\documentclass[reqno]{amsart}
\input{epsf.sty}
\usepackage{amssymb,latexsym}

\newtheorem{thm}{Theorem}[section]
\newtheorem{lemma}{Lemma}[section]

\newtheorem{prop}{Proposition}[section]
\newtheorem{defn}{Definition}[section]

\newtheorem{rem}{Remark}[section]

\makeatletter
\@addtoreset{figure}{section}
\def\thefigure{\thesection.\@arabic\c@figure}
\def\fps@figure{h, t}
\@addtoreset{table}{bsection}
\def\thetable{\thesection.\@arabic\c@table}
\def\fps@table{h, t}
\@addtoreset{equation}{section}

\makeatother

\allowdisplaybreaks
\def\intprod{\mathbin{\hbox to 6pt{%
                 \vrule height0.4pt width5pt depth0pt
                 \kern-.4pt
                 \vrule height6pt width0.4pt depth0pt\hss}}}

\begin{document}

\title[Second-order multisymplectic field theory]
{A variational approach to second-order multisymplectic field theory}
\author[S. Kouranbaeva]{Shinar Kouranbaeva}
\address{ Department of Mathematics, University of California,
Santa Cruz, CA 95064}
\email{shinar@cats.ucsc.edu}

\author[S. Shkoller]{Steve Shkoller}
\address{ Department of Mathematics, University of California,
Davis, CA 95616}
\email{shkoller@math.ucdavis.edu}

%\subjclass{Primary 58B20, 58D05; Secondary 76E99}

\date{September 1, 1998; current version March 29, 1999}
\keywords{Multisymplectic geometry, shallow water equations}

\begin{abstract}
This paper presents a geometric-variational approach to continuous and
discrete {\it second-order} field theories following the methodology of
\cite{MPS}. Staying entirely in the Lagrangian framework and letting $Y$
denote the configuration fiber bundle, we show that both the multisymplectic
structure on $J^3Y$ as well as the Noether theorem arise from the first
variation of the action function. We generalize the multisymplectic form
formula derived for first order field theories in \cite{MPS}, to the case of
second-order field theories, and we apply our theory to the 
Camassa-Holm (CH)
equation in both the continuous and discrete settings. Our discretization
produces a multisymplectic-momentum integrator, a generalization of the
Moser-Veselov rigid body  algorithm to the setting of nonlinear PDEs with
second order Lagrangians.
\end{abstract}

\maketitle

\tableofcontents

\section{Introduction}

This paper continues the development of the variational approach to
multisymplectic field theory introduced in Marsden, Patrick, and
Shkoller \cite{MPS}.  In that paper, only first-order field theories
were considered.  Herein, we shall focus on second-order field theories,
i.e., those field theories governed by Lagrangians that depend on the
spacetime location, the field, and its first and second partial derivatives.

Multisymplectic geometry and its applications to covariant field theory
and nonlinear partial differential equations (PDE) has a rich and interesting 
history that we shall not discuss in this paper; rather, we refer the reader to 
\cite{Ga74,G91a,GIM,MPS} and the references therein.  The covariant multisymplectic
approach is the field-theoretic generalization of the symplectic approach to
classical mechanics.  The configuration manifold $Q$ of classical Lagrangian 
mechanics is replaced by a fiber bundle $Y \rightarrow X$ over the $n$$+$$1$
dimensional spacetime manifold $X$, whose sections are the physical fields 
of interest; the Lagrangian phase space is $TQ$ in Lagrangian mechanics, whereas for
kth-order field theories, the role of phase space is  played by the kth-jet bundle
of $Y$, $J^kY$, thus reflecting the additional dependence of the fields on 
spatial variables.  

For a given smooth Lagrangian $L:TQ\rightarrow {\mathbb R}$, there is a 
distinguished symplectic $2$-form $\omega_L$ on $TQ$, whose Hamiltonian vector
field is the solution of the Euler-Lagrange equations of Lagrangian mechanics.  
Lagrangian field theories, on the other hand, governed by covariant Lagrangians
${\mathcal L}: J^kY \rightarrow \Lambda^{n+1}(X)$, can be completely described
by the multisymplectic $n$$+$$2$-form $\Omega_{\mathcal L}$ on $J^{2k-1}Y$, 
the field-theoretic analogue of the symplectic $2$-form $\omega_L$ of classical
mechanics. In the case that $X$ is one dimensional, $\Omega_{\mathcal L}$
reduces to the usual time-dependent $2$-form of classical nonautonomous mechanics
(see \cite{MS}).

Traditionally, the symplectic $2$-form $\omega_L$ as well as  the multisymplectic
$n$$+$$2$-form $\Omega_{\mathcal L}$ are constructed on the Lagrangian side,
using the pull-back by the Legendre transform of canonical differential forms on the 
dual or Hamiltonian side.  Recently, however, Marsden, Patrick, and Shkoller
\cite{MPS} have shown that for first order field theories wherein
${\mathcal L}: J^1Y \rightarrow \Lambda^{n+1}(X)$, 
$\Omega_{\mathcal L} = d \Theta_{\mathcal L}$ arises as the boundary term in the
first variation of the action $\int_X {\mathcal L} \circ j^1\phi$ for smooth
mappings $\phi: X \rightarrow Y$.   This method is advantageous to the traditional
approach in that 
\begin{itemize}
\item[{\bf{(i)}}] a complete geometric theory can be derived while staying entirely
on the Lagrangian side, and 
\item[{\bf{(ii)}}] multisymplectic structure can be obtained in non-standard settings
such as discrete field theory.
\end{itemize}

The purpose of this paper, is to generalize the results of \cite{MPS} to the case that
${\mathcal L}: J^2Y \rightarrow \Lambda^{n+1}(X)$.   In Section 2, we prove in
Theorem \ref{thm_main}, that a unique multisymplectic $n$$+$$2$ form arises as the
boundary term of the first variation of the action function.  We then prove in
Theorem \ref{thm_mff}, the multisymplectic form formula for second-order field
theories, a covariant generalization of the fact that in conservative mechanics,
the flow preserves the symplectic structure.  We then obtain
the covariant Noether theorem for second-order field theories, by taking the 
first variation of the action function, restricted to the space of solutions of the
covariant Euler-Lagrange equations.

In Section 3, we use our abstract geometric theory on the Camassa-Holm
(CH)
equation, a model of shallow water waves that simultaneously exhibits solitary wave 
interaction and wave-breaking.  We  show that the multisymplectic form formula produces
a new conservation law ideally suited to study wave instability, and connect our 
intrinsic theory with Bridges' theory of multisymplectic structures (see \cite{B97} and
\cite{MPS}).

Section 4 is devoted to the discretization of second-order field theories.  We
are able to use our general theory to produce numerical algorithms for nonlinear PDE
that are governed by second-order Lagrangians which naturally respect a discrete
multisymplectic form formula and a discrete Noether theorem.  Again, we demonstrate
this methodology on the  CH equation.

\section{Variational Principles for second order Classical Field
         Theory}
\subsection{Multisymplectic geometry}
In this section, we review some aspects of multisymplectic geometry
following \cite{G91a}, \cite{GIM}, \cite{MS}, and \cite{MPS}.

 Let $X$ be an orientable ($n$$+$$1$)-dimensional manifold (which in
 applications is usually spacetime) and let $\pi_{XY} :
 Y\longrightarrow X$ be a fiber bundle over $X$. Sections $\phi
 : X\longrightarrow Y$ of this \emph{covariant configuration
 bundle} will be the physical fields. The space of sections of
 $\pi_{XY}$ will be denoted $C^\infty (\pi_{XY} )$ or by
$C^\infty(Y)$. The vertical bundle $VY$
 is the subbundle $\ker T\pi_{XY}$ of $TY$, where $T\pi_{XY} $ denotes the
 tangent map of the $\pi_{XY}$.

 If $X$ has local coordinates $x^\mu $, $\mu = 1,2,\ldots , n,
 0$, adapted coordinates on $Y$ are $y^A$, $A = 1,\ldots , N$,
 along the fibers $Y_x := \pi^{-1}_{XY}(x)$, where $x\in X$ and
 $N$ is the fiber dimension of $Y$.

$J^kY$ denotes the
$k^{\rm th}$ jet bundle of $Y$, and this bundle may be defined
inductively by $J^1(\cdots (J^1Y))$. Recall that the first
jet bundle $J^1Y$ is the affine bundle over $Y$ whose fiber over $y\in Y_x$
consists of those linear mappings $\gamma : T_xX\longrightarrow T_yY $
satisfying
\[T\pi_{XY}\circ\gamma = {\rm Identity\,\, on}\, \,  T_xX \, .\]
Coordinates $(x^\mu$, $y^A)$ on $\pi_{XY}$ induce coordinates
$y^A_\mu$ on the fibers of $J^1Y$. Given $\phi \in C^\infty(Y)$,
its tangent map at $x\in X$, denoted $T_x\phi $, is an element of
$J^1Y_{\phi (x)}$.
Therefore, the map $x\longmapsto T_x\phi $ defines a section of
$J^1Y$ regarded as a bundle over $X$. This section is denoted
$j^1(\phi )$ and is called the first jet of $\phi $, or the
first prolongation of $\phi $. In coordinates, $j^1(\phi )$ is
given by
\[x^\mu \longmapsto (x^\mu , \phi^A(x^\mu ), \partial_\nu \phi^A
  (x^\mu ))\, , \]
where $\partial_\nu = \partial /\partial x^\nu$. A section of
the bundle $J^1Y\longrightarrow X$ which is the first
prolongation of the section of $Y\longrightarrow X$ is said to
be holonomic.

The first jet bundle $J^1Y$ is the appropriate configuration bundle for
first-order field theories, i.e., field theories governed by
Lagrangians which only depend on the spacetime position, the field,
and the first partial derivatives of the field.
Herein, we shall focus on second-order field theories that are governed by
Lagrangians which additionally depend on the second partial derivatives of
the fields;  thus, in second-order field theories, the Lagrangian is
defined on  $J^2Y\equiv J^1(J^1Y)$.  Let us be more specific.

\begin{defn}
The second jet bundle is the affine bundle over $J^1Y$ whose fiber at
$\gamma\in J^1Y_y$ consists of linear mappings  $s: T_xX\longrightarrow
T_{\gamma} J^1Y$ satisfying
\[T\pi_{X, J^1Y}\circ s= {\rm Identity\,\, on}\,\, T_xX \, .\]
\end{defn}
One can define the second jet prolongation of a section $\phi :X
\longrightarrow Y$, $j^2(\phi )$, as $j^1(j^1(\phi ))$, that is a map $x
\longmapsto T_xj^1(\phi )$, where $j^1(\phi )$ is regarded as a section of
$J^1Y$ over $X$. This map defines a section of $J^2Y$ regarded
as a bundle over $X$ with $j^2(\phi )(x)$ being a linear map from $T_xX$
into $T_{j^1(\phi )(x)}J^1Y$. In coordinates, $j^2(\phi )$ is given by
\[x^\mu \longmapsto (x^\mu , \phi^A(x^\mu ), \partial_{\mu_1} \phi^A
  (x^\mu ), \partial_{\mu_2}\partial_{\mu_1} \phi^A
  (x^\mu ))\, . \]
We shall also use the notation
$\phi^A_{,\mu_1\mu_2}\equiv\partial_{\mu_2}\partial_
{\mu_1} \phi^A$
for second partial derivatives. A section $\rho $ of $J^2Y\longrightarrow X$ is said to
be $2$-holonomic if $\rho = j^2(\pi_{Y,J^2Y}\circ\rho)$. Continuing
inductively, one defines the $k^{\rm th}$ jet prolongation of $\phi$,
$j^k(\phi )$, as $j^1(\cdots (j^1(\phi )))$.

Consider a second-order Lagrangian density defined as a fiber-preserving
map $\mathcal{L}: J^2Y\longrightarrow \Lambda^{n+1}(X)$, where
$\Lambda^{n+1}(X)$ is the bundle of $(n$$+$$1)$-forms on $X$. In coordinates,
we write
\[\mathcal{L}(s) = L(x^\mu , y^A, y^A_{\mu _1}, y^A_{\mu _1\mu_2})\omega\, ,\]
where $\omega = dx^1\wedge\ldots\wedge dx^n\wedge dx^0$.

For any $k$th-order Lagrangian field theory, the fundamental geometric
structure is the Cartan form $\Theta_{\mathcal L}$; this is an $(n+1)$-form
defined on $J^{2k-1}Y$ (see \cite{G91a}). For second-order field theories,
the Cartan form is defined on $J^3Y$, the covariant analogue of the phase
space in mechanics. The Euler-Lagrange equations may be written intrinsically
as
\begin{equation}\label{Euler-Lagrange}
(j^3\phi )^*(V \intprod d\Theta_{\mathcal L})=0 \ \ \forall \ V \in T(J^3Y),
\end{equation}
where $\intprod$ denotes the interior product. Traditionally, the Cartan 
form is defined  using the pull-back by the covariant Legendre transform 
of the canonical
multisymplectic $n$$+$$1$-form on the affine dual of $J^{2k-1}Y$ (see
\cite{G91a,GIM,MS}). In local coordinates, the Cartan form on $J^3Y$ is given
by
\begin{eqnarray}\label{Cartan}
\Theta_\mathcal{L}&=&\left(\frac{\partial L}{\partial
y^A_\nu}-D_\mu\left(\frac{\partial L}{\partial y^A_{\nu\mu
}}\right)\right)dy^A\wedge\omega_\nu +\frac{\partial L}{\partial
y^A_{\nu\mu}}dy^A_\nu\wedge\omega_\mu \nonumber \\ &+&
%\end{eqnarray*}
%\begin{eqnarray*}
\left( L-\frac{\partial L}{\partial
y^A_\nu}y^A_\nu + D_\mu\left(\frac{\partial L}{\partial
y^A_{\nu\mu}}\right)y^A_\nu - \frac{\partial L}{\partial
y^A_{\nu\mu}}y^A_{\nu\mu }\right)\omega \, ,
\end{eqnarray}
where $\omega _\nu = \partial_\nu \intprod\omega$ and $\omega_{\mu \nu
}=\partial_\nu \intprod\partial_\mu \intprod\omega$, etc. For a
$k^{\rm th}$-order function $f\in C^{\infty}(J^kY,{\mathbb R})$, 
the \emph{formal} partial
derivative of $f$ in the direction $x^\mu$, denoted by $D_\mu f$, is defined
by $(j^{k+1}\phi )^*(D_\mu f)=\partial_\mu (f\circ j^k\phi )$ for all $\phi
\in C^\infty(Y)$, and is a smooth function on $J^{k+1}Y$. In jet charts,
\begin{equation} \label{DERformal}
 D_\nu f = \partial_\nu f + \frac{\partial f}{\partial y^A}y^A_\nu +
\cdots + \frac{\partial f}{\partial y^A_{\mu_1\cdots\mu_k}}y^A_
{\mu_1\cdots\mu_k\nu }\ .
\end{equation}

In the next section, we shall prove that the Cartan form arises as
the boundary term in the Lagrangian variational principle.

\subsection{The variational route to the multisymplectic form.}

In this subsection we show that a multisymplectic structure is
obtained by taking the derivative of an action functional, and use this
structure to prove the multisymplectic counterpart of the fact that in
conservative mechanics, the flow of a mechanical system consists of
symplectic maps.

Let $U$ be a smooth manifold with (piecewise) smooth closed boundary.
Define the set of smooth maps
\[ \mathcal{C}^\infty = \{\phi : U\longrightarrow Y\,\, |\,\,\pi_{XY}\circ\phi :
U\longrightarrow X\,\, {\rm is\,\, an\,\,embedding}\}\, .\]

For each $\phi\in\mathcal{C}^\infty$ set $\phi_X := \pi_{XY}\circ\phi$ and
$U_X := \phi_X(U)$ so that $\phi_X : U\longrightarrow U_X$ is a
diffeomorphism. Let $\mathcal{C}$ denote the closure of
$\mathcal{C}^\infty$ in some Hilbert or Banach space norm. The
choice of topology is not crucial in this paper, and one may assume
all fields are smooth. The tangent space to the
manifold $\mathcal{C}$ at a point $\phi\in\mathcal{C}$ is given by
\[\{V\in\mathcal{C}^\infty(X,TY)\,\, |\,\,\pi_{Y,TY}\circ V=\phi\,\, {\rm and}
\,\, V_X:=T\pi_{XY}\circ V\circ\phi_X^{-1}\,\, {\rm is \,\, a \,\, vector
\,\, field \,\, on}\,\, X\}\, .\]
%Picture with U, U_X, phi, etc.

Consider $G$, the Lie group of $\pi_{XY}$-bundle automorphisms $\eta_Y :
Y\longrightarrow Y$ covering diffeomorphisms $\eta_X:X\longrightarrow X$.
%may be include diagram
\begin{defn}
The group action $\Phi : G\times\mathcal{C}\longrightarrow\mathcal{C}$
is given by
\[\Phi (\eta_Y, \phi )=\eta_Y\circ\phi\, .\]
Note that $(\eta_Y\circ\phi)_X=\eta_X\circ\phi_X$, and if $\phi\circ
\phi_X^{-1}\in C^\infty(\pi_{U_X,Y})$, then $(\eta_Y\circ\phi)\circ\phi_X^{-1}
\circ\eta_X^{-1}\in C^\infty(\pi_{\eta_X(U_X),Y})$.
\end{defn}
The fundamental problem of the classical calculus of variations is to
extremize the action functional over the space of sections of $Y
\longrightarrow X$.

\begin{defn}
The $\textbf{action functional}$ $\mathcal{S}:\mathcal{C} \rightarrow
{\mathbb R}$ is given by
\begin{equation}\label{action}
\mathcal{S}(\phi )=\int_{U_X}\mathcal{L}(j^2(\phi\circ\phi_X^{-1})) \,\,
\hbox{for all }\phi\in\mathcal{C}\, .
\end{equation}
\end{defn}

\begin{defn}\label{extremal}
$\phi \in {\mathcal C}$ is said to be \textbf{an extremum} of $\mathcal{S}$ 
if
\[\left.\frac{d}{d\lambda}\right|_{\scriptscriptstyle \lambda = 0}
\mathcal{S}(\Phi(\eta_Y^\lambda, \phi))=0\]
for all smooth paths $\lambda\mapsto\eta_Y^\lambda$ in $G$, where
for each $\lambda$, $\eta^\lambda_Y$ covers $\eta^\lambda_X$.
\end{defn}
One may associate to each $\phi^\lambda\in{\mathcal C}$, the section of
$Y$ given by
$\eta_Y^ \lambda\circ(\phi\circ\phi_X^{-1})\circ(\eta_X^\lambda)^{-1}$;
namely
$\eta_Y^ \lambda\circ(\phi\circ\phi_X^{-1})\circ(\eta_X^\lambda)^{-1}$
maps $U_X^\lambda :=\eta_X^\lambda\circ\phi_X(U)$  into $\phi^\lambda(U)$.

If we choose the curve $\phi^\lambda$ such that $\phi^0 = \phi$ and
$(d/d \lambda)|_{\lambda =0} \Phi(\eta^\lambda_Y,\phi) = V$, then we have
that
$V = (d/d \lambda)|_{\lambda =0} \phi^\lambda$ and 
$V_X=\left.\frac{d}{d\lambda}\right|_{\scriptscriptstyle \lambda = 0}
\eta_X^\lambda$. This will be used in the following:
\begin{eqnarray}\label{var1}
d\mathcal{S}_\phi \cdot V &=& \left.\frac{d}{d\lambda}\right|_{\scriptscriptstyle
\lambda = 0}\mathcal{S}(\phi^\lambda)= \left.\frac{d}{d\lambda}\right|_
{\scriptscriptstyle \lambda = 0}\int_{U_X^\lambda}\mathcal{L}(j^2(\phi^
\lambda\circ(\phi_X^\lambda)^{-1}))\nonumber \\
&=& \int_{U_X}\left.\frac{d}{d\lambda}\right|_{\scriptscriptstyle \lambda = 0}
\mathcal{L}(j^2(\phi^\lambda\circ(\phi_X^\lambda)^{-1})) + \int_{U_X}
\left.\frac{d}{d\lambda}\right|_{\scriptscriptstyle \lambda = 0}
(\eta_X^\lambda )^*\mathcal{L}(j^2(\phi\circ\phi_X^{-1})) \nonumber \\
&=& \int_{U_X}\left.\frac{d}{d\lambda}\right|_{\scriptscriptstyle \lambda = 0}
\mathcal{L}(j^2(\phi^\lambda\circ(\phi_X^\lambda)^{-1})) + \int_{U_X}
\pounds_{V_X}\mathcal{L}(j^2(\phi\circ\phi_X^{-1})) \ ,
\end{eqnarray}
where $^*$ stands for the pull-back, and $\pounds$ denotes the Lie
derivative.

Now, let $VY\subset TY$ be
the vertical subbundle; this is the bundle over $Y$ whose fibers are
given by
  \[V_yY=\{v\in T_yY\ |\ T\pi_{XY}\cdot v=0\}\ .\]
For each $\gamma\in J^1Y_y$ there exists a natural splitting $T_yY=
\mathrm{image}\ \gamma\oplus V_yY$. For example, for a vector $V\in T_\phi
\mathcal{C}$, let $\gamma = j^1(\phi\circ\phi_X^{-1})$, $V^h:= \gamma
(V_X)$, and $V^v:=V\circ\phi_X^{-1}-V^h$. Then
\[T\pi_{XY}\circ V^h = T\pi_{XY}\circ\gamma(V_X)=\mathrm{id}_{TX}(V_X)=V_X\, .\]
On the other hand, by definition, $V_X=T\pi_{XY}\circ V\circ\phi_X^{-1}$.
Therefore, $T\pi_{XY}\cdot V^v = 0$ which confirms that any vector $V\in
T_\phi\mathcal{C}$ may be decomposed into its horizontal component
\begin{equation}\label{horizontal}
V^h=T(\phi\circ\phi_X^{-1})\cdot V_X,
\end{equation}
and its vertical component
\begin{equation}
V^v=V\circ\phi_X ^{-1}-V^h \label{vertical}.
\end{equation}
\begin{rem}
Notice that $V(x) \in T_{\phi(x)}Y$ for all $x \in U$, while
$V^h$ and $V^v$ are vector fields on $U_X=\phi_X(U)$.
\end{rem}

Next, we define prolongations of automorphisms $\eta_Y$ of $Y$ and of
elements $V\in T_\phi\mathcal{C}$.
\begin{defn}
Given an automorphism $\eta_Y$ of $Y\rightarrow X$, its first prolongation
$j^1(\eta_Y): J^1Y\longrightarrow J^1Y$ is defined via
\[j^1(\eta_Y)(\gamma)= T\eta_Y\circ\gamma\circ T\eta_X^{-1}\, .\]
If $\gamma:T_xX\longrightarrow T_yY$, then $j^1(\eta_Y)(\gamma): T
_{\eta_X(x)}X\longrightarrow T_{\eta_Y(y)}Y$, with local coordinate
expression
\begin{equation}\label{1prolong}
j^1(\eta_Y)(\gamma)=\left(\eta_X^\mu, \eta_Y^A, \left(\frac{\partial\eta_Y^A}
{\partial x^\nu} + \gamma_\nu^B\frac{\partial\eta_Y^A}{\partial y^B}\right)
\frac{(\eta_X^{-1})^\nu }{\partial x^\mu}\right)\, .
\end{equation}
\end{defn}
To define the first prolongation of a vector $V\in T_\phi\mathcal{C}$,
denoted $j^1(V)$, let $\eta_Y^\lambda$ be a flow of a vector field $v$
on $Y$ with $v\circ\phi = V$.

\begin{defn}
The first prolongation $j^1(V)$ of $V$ is a vector field on $J^1Y$ given by
\[j^1(V)= \left.\frac{d}{d\lambda}\right|_{\scriptscriptstyle \lambda = 0}
j^1(\eta_Y^\lambda )\, .\]
\end{defn}

If in a coordinate chart $V=(V^\mu, V^A)$; identifying $V$ with
$V\circ\phi_X^{-1}$, we see that (\ref{1prolong})
yields the following local expression for $j^1(V)(\gamma )$:
\begin{equation}\label{Vprolong}
j^1(V)(\gamma )=\left(V^\mu , V^A, \frac{\partial V^A}{\partial x^\mu} +
\frac{\partial V^A}{\partial y^B}\gamma_\mu^B
-\gamma_\nu^A\frac{\partial V^\nu}{\partial x_\mu}\right)\, .
\end{equation}

Using induction, one can define the $k^{\rm th}$ prolongation of an automorphism
$\eta_Y$ and the $k^{\rm th}$ prolongation of a vector $V\in T_\phi
\mathcal{C}$ for all $k\geq 1$, and these will be denoted by
$j^k(\eta_Y)$ and $j^k(V)$, respectively.

\begin{defn} \label{DERvariational}
For a $k^{\rm th}$-order function $f\in C^\infty(J^kY, {\mathbb R})$, 
the variational derivative of $f$ is the function on $J^{2k}Y$ given by
\[\frac{\delta f}{\delta y^A}=\Sigma_{s=0}^k(-1)^s D_{\mu_1}\cdots
 D_{\mu_s}\left(\frac{\partial f}{\partial y^A_{\mu_1\cdots\mu_s}}\right)\ .\]
\end{defn}
In particular, for a second order function $f\in C^\infty (J^2Y,{\mathbb R})$,
the variational derivative of $f$ is the function on $J^4Y$ given by
\[
\frac{\delta f}{\delta y^A}=\frac{\partial f}{\partial y^A} - D_\nu
\left(\frac{\partial f}{\partial y^A_\nu }\right) + D_\nu D_\mu
\left(\frac{\partial f}{\partial y^A_{\nu\mu}}\right)\ .
\]

\begin{defn}
Let ${\mathcal C}^4 =\{ j^4(\phi \circ \phi_X^{-1}) | \phi \in {\mathcal C}\}$.
\end{defn}

\begin{thm}\label{thm_main}
Given a smooth Lagrangian density $\mathcal{L}:J^2Y\longrightarrow\Lambda
^{n+1}(X)$, there exist a unique $\Psi \in \Lambda^{n+2}(J^4Y)$  given by
\[ \Psi=\frac{\delta L}{\delta y^A}dy^A\wedge\omega\ ,\]
a unique map $\mathcal{D}_{EL}\mathcal{L}\in C^\infty({\mathcal C}^4,
T^\ast {\mathcal C} \otimes\Lambda^{n+1}(X))$ given by
\begin{equation}\label{DEL}
\mathcal{D}_{EL}\mathcal{L}(\phi)\cdot V = j^4(\phi\circ\phi_X^{-1})^\ast
\left(\frac{\delta L}{\delta y^A}\mathbf{i}_V (dy^A\wedge\omega )\right)\ ,
\end{equation}
and a unique differential form
$\Theta_{\mathcal{L}}\in\Lambda^{n+1}(J^3Y)$ given by
\begin{eqnarray}\label{Cartanagain}
\Theta_\mathcal{L}&=&\left(\frac{\partial L}{\partial
y^A_\nu}-D_\mu\left(\frac{\partial L}{\partial y^A_{\nu\mu
}}\right)\right)dy^A\wedge\omega_\nu +\frac{\partial L}{\partial
y^A_{\nu\mu}}dy^A_\nu\wedge\omega_\mu \nonumber \\ &+&
\left( L-\frac{\partial L}{\partial
y^A_\nu}y^A_\nu + D_\mu\left(\frac{\partial L}{\partial
y^A_{\nu\mu}}\right)y^A_\nu - \frac{\partial L}{\partial
y^A_{\nu\mu}}y^A_{\nu\mu }\right)\omega
\end{eqnarray}
such that $j^3(\phi\circ\phi_X^{-1})^\ast\Theta_{\mathcal{L}} =
\mathcal{L}\circ j^2(\phi\circ\phi_X^{-1})$ for any $\phi\in\mathcal{C}$,
and the variation of the action functional $\mathcal{S}$ is expressed by
the following formula: for any $V\in T_\phi\mathcal{C}$ and any open
subset $U_X$ of $X$ such that $\overline{U_X}\bigcap\partial X = \emptyset$,
\begin{equation}\label{var}
d\mathcal{S}_\phi \cdot V = \int_{U_X}\mathcal{D}_{EL}\mathcal{L}(\phi)\cdot V
+ \int_{\partial U_X}j^3(\phi\circ\phi_X^{-1})^\ast[j^3(V)\intprod
\Theta_{\mathcal{L}}]\, .
\end{equation}
Furthermore,
\begin{equation} \label{211}
\mathcal{D}_{EL}\mathcal{L}(\phi)\cdot V =  j^3(\phi\circ\phi_X^{-1})^\ast
[j^3(V)\intprod\Omega_{\mathcal{L}}]\ \textit{in}\  U_X \ ,
\end{equation}
where $\Omega_{\mathcal{L}} = d\Theta_{\mathcal{L}}$ is the
multisymplectic form on $J^3Y$. The variational principle (\ref{var}) yields 
the Euler-Lagrange equations (\ref{Euler-Lagrange}) on the interior of
the domain, which in coordinates are
given by
\begin{equation}\label{coordinateEL}
\frac{\partial L}{\partial y^A}(j^2(\phi\circ\phi_X^{-1}))-\frac{\partial}{\partial
x^\nu}\left(\frac{\partial L}{\partial y^A_\nu}(j^2(\phi\circ\phi_X^{-1}))\right)
+\frac{\partial^2}{\partial x^\nu\partial x^\mu}\left(\frac{\partial L}{\partial
 y^A_{\nu\mu}}(j^2(\phi\circ\phi_X^{-1}))\right) = 0\, ,
\end{equation}
while the form $\Theta_{\mathcal{L}}$ naturally arises in the boundary term
and matches the definition of the Cartan form given in (\ref{Cartan}).
\end{thm}
\begin{proof}
The proof proceeds in three steps. We begin by computing the first variation
using (\ref{var1}).
Then we show that the boundary term yields the Cartan form.
Lastly, we verify  the statements related to the interior integral.

Choose $U_X=\phi_X(U)$ small enough so that it is contained in a
coordinate chart. If in these coordinates $V=(V^\mu , V^A)$, then along $\phi\circ
\phi_X^{-1}$, the coordinate expressions for $V_X, V^h, V^v$ are written as
\begin{eqnarray} \label{vert}
V_X=V^\mu\frac{\partial}{\partial x^\mu},\quad V^h= V^\mu\frac{\partial}
{\partial x^\mu} +V^\mu\frac{\partial (\phi\circ\phi_X^{-1})^A}{\partial x^\mu}
\frac{\partial}{\partial y^A}\, ,\nonumber \\
\textrm{and}\quad V^v=(V^v)^A\frac{\partial}{\partial y^A}:= \left(
V^A - V^\mu\frac{\partial (\phi\circ\phi_X^{-1})^A}{\partial x^\mu}\right)
\frac{\partial}{\partial y^A}\, .
\end{eqnarray}

Using the Cartan formula we first compute the second term on the 
right-hand-side of (\ref{var1}):
\begin{align}\label{easyterm}
\int_{U_X}\pounds_{V_X}\mathcal{L}(j^2(\phi\circ\phi_X^{-1})) &=
\int_{U_X}\pounds_{V_X}(L\omega) \nonumber\\
&=\int_{U_X}d\mathbf{i}_{V_X}(L\omega ) + \mathbf{i}_{V_X}d(L\omega )\nonumber\\
&=\int_{\partial U_X}L\, \mathbf{i}_{V_X}\omega =\int_{\partial U_X}L\, V^\theta\omega
_\theta\, .
\end{align}

Using (\ref{vertical}),
and the local expression for the vertical vector field $V^v$, we have that
\begin{align}\label{var2}
\int_{U_X}\left.\frac{d}{d\lambda}\right|_{\scriptscriptstyle \lambda = 0}
\mathcal{L}(j^2(\phi^\lambda &\circ (\phi_X^\lambda)^{-1})) =
\int_{U_X}\left[ \frac{\partial L}{\partial y^A}(j^2(\phi\circ\phi_X^{-1}))(V^v)^A
\right. \\
&+\left. \frac{\partial L}{\partial y^A_\nu}(j^2(\phi\circ\phi_X^{-1}))(V^v)^A
_{,\nu}+ \frac{\partial L}{\partial y^A_{\nu\mu}}(j^2(\phi\circ\phi_X^{-1}))
(V^v)^A_{,\nu\mu}\right] \omega\nonumber\, .
\end{align}

In the following, we shall use $D_\nu f$ for the
formal partial derivative of a function $f$ (see (\ref{DERformal})),
and $\frac{\partial f}{\partial x^\nu}$ will denote $\frac{\partial}{\partial
x^\nu}(f\circ j^2(\phi\circ\phi_X^{-1}))$. Integrating (\ref{var2}) by
parts, we obtain that
\begin{align*}
\int_{U_X}&\left[ \frac{\partial L}{\partial y^A}-\frac{\partial}{\partial
x^\nu}\left( \frac{\partial L}{\partial y^A_\nu}\right)+\frac{\partial^2}
{\partial x_\nu\partial x_\mu}\left( \frac{\partial L}{\partial y^A_{\nu\mu}}
\right)\right](V^v)^A\omega \\
&+\int_{U_X}\left(\frac{\partial L}{\partial y^A_\nu}(V^v)^A\right)_{,\nu}
\omega +\int_{U_X}\left(\frac{\partial L}{\partial y^A_{\nu\mu}}
(V^v)^A_{,\nu}\right)_{,\mu}\omega \\
&-\int_{U_X}\left( \frac{\partial}{\partial
 x^\mu}\left(\frac{\partial L}
{\partial y^A_{\nu\mu}}\right)(V^v)^A\right)_{,\nu}\omega\, .
\end{align*}

Using the fact $f_{,\nu}\omega =d(f\omega_\nu)$, applying the Stoke's formula
$\int_Ud\alpha=\int_{\partial U}\alpha$, and  combining the last calculation
with (\ref{easyterm}), we obtain
\begin{align}\label{var3}
d\mathcal{S}_\phi \cdot V &=\int_{U_X}\left[ \frac{\partial L}{\partial y^A}
-\frac{\partial}{\partial
x^\nu}\left( \frac{\partial L}{\partial y^A_\nu}\right)+\frac{\partial^2}
{\partial x_\nu\partial x_\mu}\left( \frac{\partial L}{\partial y^A_{\nu\mu}}
\right)\right](V^v)^A\omega \\
&+\int_{\partial U_X}\left(\frac{\partial L}{\partial y^A_\nu}-\frac{\partial}
{\partial x^\mu}\left(\frac{\partial L}
{\partial y^A_{\nu\mu}}\right)\right)(V^v)^A\omega_\nu +\frac{\partial L}
{\partial y^A_{\nu\mu}}(V^v)^A_{,\nu}\omega_\mu + L\, V^\theta\omega_\theta
\, .\nonumber
\end{align}

\begin{defn}
A form $\alpha$ on $J^kY$ is contact, if $(j^k\phi )^\ast\alpha =0$
for all $\phi\in C^\infty(Y)$.
\end{defn}

\begin{lemma} \label{boundary}
For a smooth Lagrangian density $\mathcal{L}:J^2Y\longrightarrow\Lambda
^{n+1}(X)$ there exist a unique differential form
$\Theta_{\mathcal{L}}\in\Lambda^{n+1}(J^3Y)$ defined by (\ref{Cartanagain})
such that the boundary integral in (\ref{var3}) is equal to
\[\int_{\partial U_X}j^3(\phi\circ\phi_X^{-1})^\ast[j^3(V)\intprod
\Theta_{\mathcal{L}}]\, .\]
Furthermore, $\Theta_{\mathcal{L}}$ can be written as a sum of $L\omega$
and a linear combination of a system of contact forms on $J^2Y$ with
coefficients being functions on $J^3Y$.
\end{lemma}
\textit{Proof of Lemma~\ref{boundary}} Let $W=(W^\mu, W^A, W^A_\mu , W^A_{\mu\nu},
W^A_{\mu\nu\theta})$ be an arbitrary vector field on $J^3Y$, and let $\varphi:=
j^3(\phi\circ\phi_X^{-1})$, a map from $U_X$ to $J^3Y$. Then one
computes
\begin{align*}
\mathbf{i}_W(dy^A\wedge\omega_\nu )&=W^A\omega_\nu - W^\theta
dy^A\wedge\omega_{\nu\theta}\, ,\\
\mathbf{i}_W(dy^A_\nu\wedge\omega_\mu )&=W^A_\nu\omega_\mu - W^\theta
dy^A_\nu\wedge\omega_{\mu\theta}\, ,\\
\varphi^\ast\mathbf{i}_W(dy^A\wedge\omega_\nu )&=W^A\omega_\nu - W^\theta
\frac{\partial (\phi\circ\phi_X^{-1})^A}{\partial x^\mu}dx^\mu\wedge\omega_
{\nu\theta}\, .
\end{align*}
Using the formula
\begin{eqnarray} \label{formula}
dx^\mu\wedge\omega_{\nu\theta}=\left\{ \begin{array}{lll}
                                         0 &\mbox{if $\mu\neq\nu ,\theta$}\\
                                         \omega_\nu &\mbox{if $\mu =\theta$}\\
                                         -\omega_\theta &\mbox{if $\mu =\nu$,}
                                         \end{array}
                                 \right.
\end{eqnarray}
one finds that
\[\varphi^\ast\mathbf{i}_W(dy^A\wedge\omega_\nu )=W^A\omega_\nu - W^\theta
\frac{\partial (\phi\circ\phi_X^{-1})^A}{\partial x^\theta}\omega_\nu +
\frac{\partial (\phi\circ\phi_X^{-1})^A}{\partial x^\nu}W^\theta\omega_\theta\, .\]
Similarly,
\[\varphi^\ast\mathbf{i}_W(dy^A_\nu\wedge\omega_\mu )=W^A_\nu\omega_\mu -
W^\theta\frac{\partial^2(\phi\circ\phi_X^{-1})^A}{\partial x^\theta\partial
x^\nu}\omega_\mu + \frac{\partial^2(\phi\circ\phi_X^{-1})^A}{\partial x^\mu
\partial x^\nu}W^\theta\omega_\theta\, .\]
Thus if we let $W=j^3(V)$,
use (\ref{Vprolong}), and recall the local expression
(\ref{vert}) for $(V^v)^A$, we obtain that
\begin{align*}
\varphi^\ast\mathbf{i}_{j^3(V)}(dy^A\wedge\omega_\nu )& = (V^v)^A\omega_\nu +
\frac{\partial (\phi\circ\phi_X^{-1})^A}{\partial x^\nu}V^\theta\omega_\theta\, ,\\
\varphi^\ast\mathbf{i}_{j^3(V)}(dy^A_\nu\wedge\omega_\mu )
&=(V^v)^A_{,\nu}\,
\omega_\mu +\frac{\partial^2(\phi\circ\phi_X^{-1})^A}{\partial x^\mu
\partial x^\nu}V^\theta\omega_\theta\, .
\end{align*}
Next, observe that $V^\theta\omega_\theta = \mathbf{i}_V\omega$. Also,
\[\frac
{\partial (\phi\circ\phi_X^{-1})^A}{\partial x^\nu}=j^3(\phi\circ\phi_X^{-1})^\ast
 y^A_\nu ,\ \mbox{and}\ \frac{\partial^2(\phi\circ\phi_X^{-1})^A}{\partial x^\mu
\partial x^\nu}=j^3(\phi\circ\phi_X^{-1})^\ast y^A_{\nu\mu}\, .\]
These observations together with the previous identities imply the following
important formulas
\begin{equation}\label{VV}
\begin{array}{c}
j^3(\phi\circ\phi_X^{-1})^\ast [j^3(V)\intprod (dy^A_\nu\wedge\omega_\mu -
y^A_{\nu\mu}\omega)]= (V^v)^A_{,\nu}\,\omega_\mu\, ,\\
j^3(\phi\circ\phi_X^{-1})^\ast [j^3(V)\intprod (dy^A\wedge\omega_\nu -
y^A_{\nu}\omega)]= (V^v)^A\omega_\nu\, .
\end{array}
\end{equation}
Substituting these formulas into the boundary integral of the variational
principle (\ref{var3}), we obtain that
\begin{align*}
\int_{\partial U_X}\left(\frac{\partial L}{\partial y^A_\nu}-\frac{\partial}
{\partial x^\mu}\left(\frac{\partial L}
{\partial y^A_{\nu\mu}}\right) \right)&(V^v)^A\omega_\nu +\frac{\partial L}
{\partial y^A_{\nu\mu}}(V^v)^A_{,\nu}\omega_\mu + L\, V^\theta\omega_\theta \\
=\int_{\partial U_X}j^3(\phi\circ\phi_X^{-1})^\ast [j^3(V)\intprod &
\left\{
\left(
\frac{\partial L}{\partial
y^A_\nu}-D_\mu
\left(\frac{\partial L}{\partial y^A_{\nu\mu }}
\right)
\right)(dy^A\wedge\omega_\nu - y^A_\nu\omega)
\right. \\
&+
\left.
\left.
\frac{\partial L}{\partial y^A_{\nu\mu}}(dy^A_\nu\wedge\omega_\mu -
y^A_{\nu\mu }\omega) + L\omega
\right\}
\right] \\
=\int_{\partial U_X}j^3(\phi\circ\phi_X^{-1})^\ast [j^3(V)\intprod &
\left\{
\left(
\frac{\partial L}{\partial
y^A_\nu}-D_\mu
\left(\frac{\partial L}{\partial y^A_{\nu\mu }}
\right)
\right)dy^A\wedge\omega_\nu +\frac{\partial L}{\partial
y^A_{\nu\mu}}dy^A_\nu\wedge\omega_\mu
\right. \\
&+
\left.
\left.
\left( L-\frac{\partial L}{\partial
y^A_\nu}y^A_\nu + D_\mu
\left(\frac{\partial L}{\partial
y^A_{\nu\mu}}
\right)y^A_\nu - \frac{\partial L}{\partial
y^A_{\nu\mu}}y^A_{\nu\mu }
\right)\omega
\right\}
\right] \\
=\int_{\partial U_X}j^3(\phi\circ\phi_X^{-1})^\ast [j^3(V)\intprod
\Theta_\mathcal{L}] \, .
\end{align*}
This proves the existence of a unique differential form $\Theta_\mathcal{L}$
and demonstrates how this form naturally arises in the boundary integral of
the variational principle. Integration by parts yields the boundary integral
with terms that involve partial derivatives of $(V^v)^A$ of all orders up to
$k-1$ (in our case $k=2$). Equation (\ref{VV}) show that each partial
derivative of $(V^v)^A$ has an associated $(n+1)$-form on $J^2Y$, and
substitution of these forms yields a unique differential $(n+1)$-form as
desired. Since $L$ and its partial derivatives are functions on $J^2Y$, then
by (\ref{DERformal}), $D_\mu\left(\frac{\partial L}{\partial
y^A_{\nu\mu}}\right)$ is a function on $J^3Y$, and therefore
$\Theta_\mathcal{L}$ is a $(n+1)$-form on $J^3Y$.

It is easy to show that
\begin{align*}
j^k(\phi\circ\phi_X^{-1})^\ast (dy^A\wedge\omega_\nu - y^A_\nu\omega)=0 \\
j^k(\phi\circ\phi_X^{-1})^\ast (dy^A_\nu\wedge\omega_\mu -
y^A_{\nu\mu }\omega)=0
\end{align*}
for all integers $k\geq 2$ and for all $\phi\in \mathcal{C}$.
Therefore, $dy^A\wedge\omega_\nu - y^A_\nu\omega$ and
$dy^A_\nu\wedge\omega_\mu - y^A_{\nu\mu }\omega$ are contact forms on
$J^2Y$. Hence the last statement of the lemma follows. $\square$

A simple computation then verifies that
$\Theta_\mathcal{L}$ is the Cartan form so that
\[j^3(\phi\circ\phi_X^{-1})^\ast\Theta_{\mathcal{L}} =
   \mathcal{L}\circ j^2(\phi\circ\phi_X^{-1})\, .\]

Next, consider the interior integral of the variational
principle (\ref{var3}). Since $j^k(\phi\circ\phi_X^{-1})^\ast
\mathbf{i}_{j^k(V)}(dy^A\wedge\omega) = (V^v)^A\omega$ for all integers
$k\geq 1$, we obtain that
\begin{align} \label{interior}
&\int_{U_X}\left[ \frac{\partial L}{\partial y^A}
-\frac{\partial}{\partial
x^\nu}\left( \frac{\partial L}{\partial y^A_\nu}\right)+\frac{\partial^2}
{\partial x_\nu\partial x_\mu}\left( \frac{\partial L}{\partial y^A_{\nu\mu}}
\right)\right](V^v)^A\omega \nonumber \\
&=\int_{U_X}j^4(\phi\circ\phi_X^{-1})^\ast\mathbf{i}_{j^4(V)}
\left[ \frac{\partial L}{\partial y^A}
-D_\nu\left( \frac{\partial L}{\partial y^A_\nu}\right)+D_\nu D_\mu\left(
\frac{\partial L}{\partial y^A_{\nu\mu}}\right)\right]dy^A\wedge\omega
\nonumber \\
&=\int_{U_X}j^4(\phi\circ\phi_X^{-1})^\ast\mathbf{i}_{j^4(V)}
\left(\frac{\delta L}{\delta y^A}dy^A\wedge\omega\right)\, ,
\end{align}
where $\frac{\delta L}{\delta y^A}$ is the variational derivative of $L$ in the
direction $y^A$ (see Definition ~\ref{DERvariational}).
Since $L$ is a function of second order by hypothesis,
then its variational derivative is a function on $J^4Y$. Therefore the
form $\Psi\equiv\frac{\delta L}{\delta y^A}dy^A\wedge\omega$ is an $(n+2)$-form on
$J^4Y$. Moreover, the integrand in (\ref{interior}) written as
 $j^4(\phi\circ\phi_X^{-1})^\ast
\left(\frac{\delta L}{\delta y^A}\mathbf{i}_V (dy^A\wedge\omega)\right)$
defines a unique smooth section $\mathcal{D}_{EL}\mathcal{L}\in
C^\infty(\mathcal{C}^4, T^\ast {\mathcal C} \otimes\Lambda^{n+1}(X))$ as 
desired in the statement of the theorem.

Now we shall prove the following
\begin{lemma} \label{criticalpoint}
The forms $\Omega_{\mathcal{L}}=d\Theta_{\mathcal{L}}$ and $\Psi =
\frac{\delta L}{\delta y^A}dy^A\wedge\omega$ satisfy the following
relationship:
\begin{equation} \label{PsiOmega}
j^4(\phi\circ\phi_X^{-1})^\ast\mathbf{i}_{j^4(V)}\Psi =
j^3(\phi\circ\phi_X^{-1})^\ast\mathbf{i}_{j^3(V)}\Omega_{\mathcal{L}}
\end{equation}
for all $\phi\in\mathcal{C}$ and all vectors $V\in T\mathcal{C}$.

Furthermore, a necessary condition for $\phi\in \mathcal{C}$ to be an
extremum of the action functional $\mathcal{S}$ is that
\begin{equation} \label{Omegacondition}
j^3(\phi\circ\phi_X^{-1})^\ast\mathbf{i}_W\Omega_{\mathcal{L}} =0
\end{equation}
for all vector fields $W$ on $J^3Y$, which is equivalent
to
\begin{equation} \label{Psicondition}
j^4(\phi\circ\phi_X^{-1})^\ast\mathbf{i}_V\Psi =0
\end{equation}
for all vector fields $V$ on $J^4Y$.
\end{lemma}
\textit{Proof of Lemma~\ref{criticalpoint}}.
The proof will involve some lengthy computations that we partially present
below. To compute $\Omega_{\mathcal{L}}$, let us write
$\Theta_{\mathcal{L}}$ as
\begin{align*}
\Theta_{\mathcal{L}} &=
\left(
\frac{\partial L}{\partial
y^A_\nu}-D_\mu
\left(\frac{\partial L}{\partial y^A_{\nu\mu }}
\right)
\right)(dy^A\wedge\omega_\nu - y^A_\nu\omega) \\
&+ \frac{\partial L}{\partial y^A_{\nu\mu}}(dy^A_\nu\wedge\omega_\mu -
y^A_{\nu\mu }\omega) + L\omega .
\end{align*}
Then, for $W\in TJ^3Y$, we obtain
\begin{align*}
\mathbf{i}_W\Omega_\mathcal{L} &= W
\left[
\frac{\partial L}{\partial
y^A_\nu}-D_\mu
\left(\frac{\partial L}{\partial y^A_{\nu\mu }}
\right)
\right](dy^A\wedge\omega_\nu - y^A_\nu\omega) \\ &- d
\left(
\frac{\partial L}{\partial y^A_\nu}-D_\mu
\left(\frac{\partial L}{\partial y^A_{\nu\mu }}
\right)
\right)\wedge (W^A\omega_\nu - W^\theta dy^A\wedge\omega_{\nu\theta }
- y^A_\nu W^\theta\omega_\theta ) \\ &+ W
\left[
\frac{\partial L}{\partial y^A_{\nu\mu}}
\right]
(dy^A_\nu\wedge\omega_\mu - y^A_{\nu\mu }\omega) \\ &- d
\left(
\frac{\partial L}{\partial y^A_{\nu\mu}}
\right)
\wedge (W^A_\nu\omega_\mu - W^\theta dy^A_\nu\wedge\omega_{\mu\theta }
- y^A_{\nu\mu } W^\theta\omega_\theta ) \\ &+ D_\mu
\left(
\frac{\partial L}{\partial y^A_{\nu\mu}}
\right)
(W^A_\nu\omega - W^\theta dy^A_\nu\wedge\omega_\theta ) +
\frac{\partial L}{\partial y^A}(W^A\omega - W^\theta dy^A\wedge\omega_\theta
)\ .
\end{align*}
The last step is to pull back $\mathbf{i}_W\Omega_\mathcal{L}$ by
$\varphi = j^3(\phi\circ\phi_X^{-1})$; this eliminates the terms with the
contact forms. In addition, using the fact that the pull back commutes
with the exterior derivative, and applying formulas such as
(\ref{formula}), we obtain that
\begin{eqnarray*}
\varphi^\ast\mathbf{i}_W\Omega_\mathcal{L} %\begin{array}[t]{ll}
&=& W^A
\left\{
\frac{\partial L}{\partial y^A}\omega - d
\left(
\frac{\partial L}{\partial y^A_\nu} - \frac{\partial}{\partial x^\mu}
\left(
\frac{\partial L}{\partial y^A_{\nu\mu}}
\right)
\right)
\wedge \omega_\nu
\right\} \\
&+& W^\theta  \begin{array}[t]{l}
            \left\{ d
            \left(
            \frac{\partial L}{\partial y^A_\nu} - \frac{\partial}
            {\partial x^\mu}
            \left(
            \frac{\partial L}{\partial y^A_{\nu\mu}}
            \right)
            \right)
            \wedge
            \left(
            \frac{\partial(\phi\circ\phi_X^{-1})^A}{\partial
            x^\theta }\omega_\nu - \frac{\partial(\phi\circ\phi_X^{-1})^A}
            {\partial x^\nu }\omega_\theta
            \right.
            \right. \\
            \left. + \frac{\partial(\phi\circ\phi_X^{-1})^A}
            {\partial x^\nu }\omega_\theta
            \right) + d
            \left(
            \frac{\partial L}{\partial y^A_{\nu\mu}}
            \right)
            \wedge
            \left(
            \frac{\partial^2(\phi\circ\phi_X^{-1})^A}{\partial
            x^\theta \partial x^\nu}\omega_\mu - \frac{\partial^2
            (\phi\circ\phi_X^{-1})^A}
            {\partial x^\mu \partial x^\nu }\omega_\theta
            \right. \\
            \left.
             + \frac{\partial^2(\phi\circ\phi_X^{-1})^A}
            {\partial x^\mu \partial x^\nu }\omega_\theta
            \right)
            \left.
            - \frac{\partial}{\partial x^\mu}
            \left(
            \frac{\partial L}{\partial y^A_{\nu\mu}}
            \right)
            \frac{\partial^2(\phi\circ\phi_X^{-1})^A}{\partial
            x^\theta \partial x^\nu}\omega - \frac{\partial L}{\partial y^A}
            \frac{\partial(\phi\circ\phi_X^{-1})^A}{\partial
            x^\theta }\omega
            \right\}
            \end{array}\\
&+& W^A_\nu
\left\{
\frac{\partial}{\partial x^\mu}
\left(
\frac{\partial L}{\partial y^A_{\nu\mu}}
\right)
\omega - d
\left(
\frac{\partial L}{\partial y^A_{\nu\mu}}
\right)
\wedge\omega_\mu
\right\} \ .
%     \end{array}
\end{eqnarray*}
Some cancellation and further rearrangement yields
\begin{align*}
\varphi^\ast\mathbf{i}_W\Omega_\mathcal{L}
&= W^A
\left(
\frac{\partial L}{\partial y^A} - \frac{\partial}{\partial x^\nu}
\left(
\frac{\partial L}{\partial y^A_\nu}
\right)
+ \frac{\partial^2}{\partial x_\nu \partial x_\mu}
\left(
\frac{\partial L}{\partial y^A_{\nu\mu}}
\right)
\right)\omega  \\
&-W^\theta \frac{\partial(\phi\circ\phi_X^{-1})^A}{\partial x^\theta}
\left(
\frac{\partial L}{\partial y^A} - \frac{\partial}{\partial x^\nu}
\left(
\frac{\partial L}{\partial y^A_\nu}
\right)
+ \frac{\partial^2}{\partial x_\nu \partial x_\mu}
\left(
\frac{\partial L}{\partial y^A_{\nu\mu}}
\right)
\right)\omega \ .
\end{align*}
Letting $W=j^3(V)$, we have that
\[
\varphi^\ast\mathbf{i}_{j^3(V)}\Omega_\mathcal{L} = (V^v)^A
\left(
\frac{\partial L}{\partial y^A} - \frac{\partial}{\partial x^\nu}
\left(
\frac{\partial L}{\partial y^A_\nu}
\right)
+ \frac{\partial^2}{\partial x_\nu \partial x_\mu}
\left(
\frac{\partial L}{\partial y^A_{\nu\mu}}
\right)
\right)\omega \ ,
\]
where the right hand side equals
$j^4(\phi\circ\phi_X^{-1})^\ast\mathbf{i}_{j^4(V)}\Psi$ by (~\ref{interior}).
Hence, the relation (\ref{PsiOmega}) is proved.

A necessary condition for $\phi\in\mathcal{C}$ to be an extremum of the
action functional $\mathcal{S}$ is that the interior integral in
(\ref{var3}) vanish for all vectors $V\in T\mathcal{C}$. From the
calculation above, one may readily see that it is equivalent to the
condition (\ref{Omegacondition}).

Now if we let $V$ be a vector field on $J^4Y$, then
\begin{align*}
\mathbf{i}_V\Psi
&= \mathbf{i}_V
\left(
\frac{\delta L}{\delta y^A} dy^A\wedge\omega
\right) \\
&= \frac{\delta L}{\delta y^A}\mathbf{i}_V(dy^A\wedge\omega) \\
&= \frac{\delta L}{\delta y^A}(\mathbf{i}_V(dy^A)\wedge\omega -
dy^A\wedge\mathbf{i}_V\omega ) \\
&= \frac{\delta L}{\delta y^A}(V^A\omega - V^\theta dy^A\wedge\omega_\theta
)\ .
\end{align*}
Hence,
\begin{align*}
&j^4(\phi\circ\phi_X^{-1})^\ast\mathbf{i}_V\Psi =\\
&
\left(
\frac{\partial L}{\partial y^A} - \frac{\partial}{\partial x^\nu}
\left(
\frac{\partial L}{\partial y^A_\nu}
\right)
+ \frac{\partial^2}{\partial x_\nu \partial x_\mu}
\left(
\frac{\partial L}{\partial y^A_{\nu\mu}}
\right)
\right)
\left(
V^A -V^\theta \frac{\partial(\phi\circ\phi_X^{-1})^A}{\partial x^\theta}
\right) \omega\ .
\end{align*}
Thus, the condition
\[
j^4(\phi\circ\phi_X^{-1})^\ast\mathbf{i}_V\Psi =0
\]
for all vector fields $V$ on $J^4Y$ is equivalent to the condition
(\ref{Omegacondition}). This completes the proof of the lemma. $\square$

Lemma~\ref{criticalpoint} contains two equivalent conditions for
$\phi\in\mathcal{C}$ to be extremal. Both conditions yield the same
coordinate expression of the Euler-Lagrange equations given by
\[
\frac{\partial L}{\partial y^A}(j^2(\phi\circ\phi_X^{-1}))-\frac{\partial}{\partial
x^\nu}\left(\frac{\partial L}{\partial y^A_\nu}(j^2(\phi\circ\phi_X^{-1}))\right)
+\frac{\partial^2}{\partial x^\nu\partial x^\mu}\left(\frac{\partial L}{\partial
 y^A_{\nu\mu}}(j^2(\phi\circ\phi_X^{-1}))\right) = 0\ ,
\]
which is the final statement of the theorem.
\end{proof}
\begin{rem}
As one may see the proof we have presented can be generalized to
Lagrangian densities on $J^kY$. One has to modify the labeling of
variables to reflect the general case. For example,
\[
(V^v)^A_{,\mu_1\ldots\mu_l}\omega_\theta =
\varphi^\ast[j^3(V)\intprod (dy^A_{\mu_1\ldots\mu_l}\wedge\omega_\theta
- y^A_{\mu_1\ldots\mu_l\theta }\omega)]\ ,
\]
where $0\leq l\leq (k-1)$. Then the Cartan form shall arise in the
boundary integral as a linear combination of the forms above.
\end{rem}

We shall call critical points $\phi$ of $\mathcal{S}$, solutions of
the Euler-Lagrange equations.
\begin{defn}
We let
\begin{equation}
\mathcal{P} = \{ \phi \in \mathcal{C}\quad |\quad
j^3(\phi\circ\phi_X^{-1})^\ast\mathbf{i}_W\Omega_\mathcal{L} = 0
\quad
\mbox{for all vector fields   }
\ W \mbox{on   }\ J^3Y \}
\end{equation}
denote the space of solutions of the Euler-Lagrange equations.
\end{defn}
We are now ready to prove the multisymplectic form formula, a
covariant generalization of the symplectic flow theorem to second-order
field theories.
\footnote{For first-order field theories, this is Theorem 4.1 in \cite{MPS}.}

\subsection{The multisymplectic form formula}
 If $\phi^\lambda$ is a smooth curve of solutions of the Euler-Lagrange
 equations in $\mathcal{P}$ (when such solutions exist), then
 differentiating with respect to $\lambda$ at $\lambda = 0$ will give a
 tangent vector $V$ to the curve at $\phi =\phi ^0$. By differentiating
 $\left.\frac{d}{d\lambda}\right|_{\scriptscriptstyle \lambda =
 0}j^3(\phi^\lambda\circ (\phi_X^\lambda)^{-1})^\ast [W\intprod\Omega_
 {\mathcal{L}}]=0$, we obtain
 \[j^3(\phi\circ\phi_X^{-1})^\ast \pounds_{j^3(V)}[W\intprod\Omega_
 {\mathcal{L}}]=0\]
 for all vector fields $W$ on $J^3Y$. Therefore if $\mathcal{P}$ is a
 submanifold of $\mathcal{C}$, then for any $\phi\in\mathcal{P}$ we may
 identify $T_\phi\mathcal{P}$ with the set of vectors $V$ that satisfy
 the above condition. However, we do not require $\mathcal{P}$ to be a
 submanifold.
 \begin{defn}
For any $\phi\in \mathcal{P}$,
\begin{equation} \label{firstvar}
\mathcal{F}=\{V\in T_\phi\mathcal{C}\quad |\quad j^3(\phi\circ\phi_X^{-1})
^\ast \pounds_{j^3(V)}[W\intprod\Omega_{\mathcal{L}}]=0 \quad
\mbox{for all vector fields }V\ \mbox{on}\ J^3Y\}
\end{equation}
defines a set of solutions of the first variation equations of the
Euler-Lagrange equations.
 \end{defn}
\begin{thm}\textbf{(Multisymplectic form formula)}\label{thm_mff}
If $\phi\in\mathcal{P}$, then for all $V$ and $W$ in $\mathcal{F}$,
\begin{equation}\label{mff}
\int_{\partial U_X} j^3(\phi\circ\phi_X^{-1})^\ast [j^3(V)\intprod
j^3(W)\intprod\Omega_{\mathcal{L}}]=0.
\end{equation}
\end{thm}
\begin{proof}
We follow Theorem 4.1 in \cite{MPS}.
Define the $1$-forms $\alpha_1$ and $\alpha_2$ on $\mathcal{C}$ by
\[
\alpha_1(\phi )\cdot V:=\int_{U_X}j^3(\phi\circ\phi_X^{-1})^\ast [j^3(V)
\intprod \Omega_{\mathcal{L}}]
\]
and
\[
\alpha_2(\phi )\cdot V:=\int_{\partial U_X}j^3(\phi\circ\phi_X^{-1})^\ast
[j^3(V)\intprod \Theta_{\mathcal{L}}]\ ,
\]
so that by (\ref{var}) and (\ref{211}),
\begin{equation}\label{dS}
d\mathcal{S}_\phi\cdot V = \alpha_1(\phi )\cdot V + \alpha_2(\phi )\cdot V
\quad \hbox{for all}\quad V\in T_{\phi}\mathcal{C}\ .
\end{equation}
Furthermore,
\[
d^2\mathcal{S}(\phi )(V, W)=d\alpha_1(\phi )(V, W) + d\alpha_2(\phi )
(V, W)\quad\hbox{for all}\quad V, W \in T_\phi\mathcal{C}\ .
\]
Since $d^2\mathcal{S}=0$, we have that
\begin{equation}\label{star}
d\alpha_1(\phi )(V, W) + d\alpha_2(\phi )(V, W)=0
\quad\hbox{for all}\quad V, W \in T_\phi\mathcal{C}\ .
\end{equation}
Given vectors $V, W\in T_\phi\mathcal{C}$ we may extend them to vector
fields $\mathcal{V}, \mathcal{W}$ on $\mathcal{C}$ by fixing vector fields
 $v, w$ on $Y$ such that $V=v\circ\phi$ and $W=w\circ\phi$, and letting
 $\mathcal{V}(\rho )=v\circ\rho$ and $\mathcal{W}(\rho )=w\circ\rho$. If
 $\eta_Y^\lambda$ covering $\eta_X^\lambda$ is the flow of $v$, then
 $\Phi(\eta_Y^\lambda , \rho )$ is the flow of $\mathcal{V}$. Notice that
 $\mathcal{V}(\phi )=V$ and $\mathcal{W}(\phi )=W$, hence the equation
(\ref{star}) becomes
\[
d\alpha_1(\mathcal{V}, \mathcal{W})(\phi ) + d\alpha_2(\mathcal{V},
\mathcal{W})(\phi )=0\ .
\]
 Recall that for any $1$-form $\alpha $ on $\mathcal{C}$ and vector fields
 $\mathcal{V}, \mathcal{W}$ on $\mathcal{C}$,
\begin{equation}\label{bracket}
d\alpha (\mathcal{V}, \mathcal{W})= \mathcal{V}[\alpha (\mathcal{W})] -
\mathcal{W}[\alpha (\mathcal{V})]-\alpha([\mathcal{V}, \mathcal{W}])\ .
\end{equation}
Also recall that for a vector field $\mathcal{V}$ on $\mathcal{C}$ and
a function $f$ on $\mathcal{C}$, $\mathcal{V}[f]=df\cdot\mathcal{V}$. We
now use the latter and (\ref{bracket}) on $\alpha_2$.  We have that
\begin{align}\label{doublestar}
d\alpha_2(\mathcal{V}, \mathcal{W})(\phi ) & = \mathcal{V}[\alpha_2
(\mathcal{W})](\phi ) - \mathcal{W}[\alpha_2(\mathcal{V})](\phi )-
\alpha_2([\mathcal{V}, \mathcal{W}])(\phi ) \nonumber \\
& =[d(\alpha_2(\mathcal{W}))\cdot\mathcal{V}](\phi ) - [d(\alpha_2
(\mathcal{V}))\cdot\mathcal{W}](\phi )- \alpha_2(\phi )\cdot [V, W] \nonumber \\
& =d(\alpha_2(\mathcal{W}))(\phi )\cdot V - d(\alpha_2
(\mathcal{V}))(\phi )\cdot W - \alpha_2(\phi )\cdot [V, W]\ .
\end{align}
Similarly,
\begin{align}\label{triplestar}
d\alpha_1(\phi )(V, W)=d(\alpha_1(\mathcal{W}))(\phi )\cdot V - d(\alpha_1
(\mathcal{V}))(\phi )\cdot W - \alpha_1(\phi )\cdot [V, W]\ .
\end{align}
Let $\phi\in\mathcal{P}$ and $\phi^\lambda =\eta_Y^\lambda\circ\phi$ be a
curve in $\mathcal{C}$ through $\phi$ such that
\[
V=\left.\frac{d}{d\lambda}\right|_{\scriptscriptstyle \lambda
=0}\phi^\lambda\quad\hbox{and}\quad V\in\mathcal{F}\ .
\]
Now we restrict $V, W$ to
$\mathcal{F}$. We shall give a detailed computation of the first term on the right hand
side of (\ref{doublestar}).  We have that
\begin{align*}
d(\alpha_2(\mathcal{W}))(\phi )\cdot V & =
\left.\frac{d}{d\lambda}\right|_{\scriptscriptstyle \lambda =0}
(\alpha_2(\mathcal{W}))(\phi^\lambda )\\
& =\left.\frac{d}{d\lambda}\right|_{\scriptscriptstyle \lambda =0}
\alpha_2(\phi^\lambda )\cdot (w\circ\phi^\lambda )\\
& =\left.\frac{d}{d\lambda}\right|_{\scriptscriptstyle \lambda =0}
\int_{\partial (\eta_X^\lambda (U_X))}j^3(\phi^\lambda\circ(\phi_X
^\lambda )^{-1})^\ast [j^3(w\circ\phi^\lambda )\intprod
\Theta_{\mathcal{L}}]\ \\
& = \left.\frac{d}{d\lambda}\right|_{\scriptscriptstyle \lambda =0}
\int_{\partial U_X}j^3(\phi\circ\phi_X^{-1})^\ast j^3(\eta_Y^\lambda )^\ast
[j^3(W)\intprod \Theta_{\mathcal{L}}] \\
& =\int_{\partial U_X}j^3(\phi\circ\phi_X^{-1})^\ast \pounds_{j^3(V)}
(j^3(W)\intprod \Theta_{\mathcal{L}})\\
& =\int_{\partial U_X}j^3(\phi\circ\phi_X^{-1})^\ast d[j^3(V)\intprod
j^3(W)\intprod \Theta_{\mathcal{L}}]\\
& +\int_{\partial U_X}j^3(\phi\circ\phi_X^{-1})^\ast [j^3(V)\intprod
d(j^3(W)\intprod \Theta_{\mathcal{L}})]\ ,
\end{align*}
where the last equality was obtained using Cartan's formula. We have also
used the fact that $W^\lambda =w\circ\phi^\lambda$ and $W=w\circ\phi$ have
the same $k$-th prolongation. Furthermore, using Stoke's theorem, noting
that $\partial\partial U_X$ is empty, and applying Cartan's formula once
again to $d(j^3(W)\intprod \Theta_{\mathcal{L}})$, we obtain that
\begin{align}\label{one}
d(\alpha_2(\mathcal{W}))(\phi )\cdot V  & =
\int_{\partial U_X}j^3(\phi\circ\phi_X^{-1})^\ast [j^3(V)\intprod
\pounds_{j^3(W)}\Theta_{\mathcal{L}}] \nonumber \\
& -\int_{\partial U_X}j^3(\phi\circ\phi_X^{-1})^\ast [j^3(V)\intprod
j^3(W)\intprod \Omega_{\mathcal{L}}]\ .
\end{align}
Similarly,
\begin{align}\label{two}
d(\alpha_2(\mathcal{V}))(\phi )\cdot W  & =
\int_{\partial U_X}j^3(\phi\circ\phi_X^{-1})^\ast [j^3(W)\intprod
\pounds_{j^3(V)}\Theta_{\mathcal{L}}] \nonumber \\
& -\int_{\partial U_X}j^3(\phi\circ\phi_X^{-1})^\ast [j^3(W)\intprod
j^3(V)\intprod \Omega_{\mathcal{L}}]\ .
\end{align}
Now, $j^3([V, W])= [j^3(V), j^3(W)]$; hence,
\[
\alpha_2(\phi )\cdot[V, W]= \int_{\partial U_X}j^3(\phi\circ\phi_X^{-1})
^\ast ([j^3(V), j^3(W)]\intprod\Theta_{\mathcal{L}})\ .
\]
Recall that for a differential form $\alpha$ on a manifold $M$ and for vector
fields $X, Y$ on $M$,
\[
\mathbf{i}_{[X, Y]}\alpha = \pounds_X\mathbf{i}_Y\alpha - \mathbf{i}_Y
\pounds_X\alpha\ .
\]
Therefore,
\begin{align}\label{three}
\alpha_2(\phi )\cdot[V, W]& = \int_{\partial U_X}j^3(\phi\circ\phi_X^{-1})
^\ast [\pounds_{j^3(V)}(j^3(W)\intprod\Theta_{\mathcal{L}})-
j^3(W)\intprod\pounds_{j^3(V)}\Theta_{\mathcal{L}}]\nonumber \\
 =\int_{\partial U_X}j^3(\phi\circ & \phi_X^{-1})^\ast [j^3(V)\intprod
\pounds_{j^3(W)}\Theta_{\mathcal{L}}-j^3(V)\intprod j^3(W)
\intprod\Omega_{\mathcal{L}}-j^3(W)\intprod\pounds_{j^3(V)}\Theta_{\mathcal{L}}],
\end{align}
where we have again used Stoke's theorem and Cartan's formula twice.
Substituting (\ref{one}), (\ref{two}), and (\ref{three}) into
(\ref{doublestar}), we obtain that
\begin{equation}\label{four}
d\alpha_2(\phi )(V, W)=\int_{\partial U_X}j^3(\phi\circ\phi_X^{-1})^\ast
[j^3(W)\intprod j^3(V)\intprod\Omega_{\mathcal{L}}]\ .
\end{equation}
We now compute (~\ref{triplestar}). Similar computations as above yield
\[
d(\alpha_1(\mathcal{W}))(\phi )\cdot V=\int_{U_X}j^3(\phi\circ\phi_X^{-1})
^\ast\pounds_{j^3(V)}(j^3(W)\intprod\Omega_{\mathcal{L}})
\]
which vanishes for all $\phi\in\mathcal{P}$ and $V\in\mathcal{F}$.
Similarly, $d(\alpha_1(\mathcal{V}))(\phi )\cdot W=0$ for all $\phi\in
\mathcal{P}$ and $W\in\mathcal{F}$. Finally, $\alpha_1(\phi )=0$ for all
$\phi\in\mathcal{P}$. Therefore, the equation (~\ref{triplestar})
vanishes for all $\phi\in\mathcal{P}$ and $V, W \in\mathcal{F}$. Using the
latter and (\ref{four}), equation (\ref{star}) becomes
\[
\int_{\partial U_X} j^3(\phi\circ\phi_X^{-1})^\ast [j^3(W)\intprod
j^3(V)\intprod\Omega_{\mathcal{L}}]=0
\]
for all $\phi\in\mathcal{P}$ and all $V, W \in \mathcal{F}$, as desired.
 \end{proof}

\subsection{The Noether Theorem}
Suppose that $\mathcal{S}$ is invariant under the action $\Phi(g, \phi)$ of a
Lie group $G$ on $\mathcal{C}$. This implies that for each $g\in G$, $\Phi(g,
\phi)\in \mathcal{P}$ whenever $\phi\in \mathcal{P}$. We restrict the action
to elements of $\mathcal{P}$. For each element $\xi$ of the
Lie algebra ${\mathfrak g}$ of $G$, let $\xi_\mathcal{C}$ be the corresponding
infinitesimal
generator on $\mathcal{C}$ restricted to elements of $\mathcal{P}$. By the
invariance of $\mathcal{S}$,
\[
\mathcal{S}(\Phi(\mathrm{exp}(t\xi ), \phi ))=\mathcal{S}(\phi )\quad
\mbox{for all}\quad t\, .
\]
Differentiating with respect to $t$ at $t=0$, and using the fundamental
property of the Cartan form that $\mathcal{L}\circ
j^2(\phi\circ\phi_X^{-1})=j^3(\phi\circ\phi_X^{-1})^\ast\Theta_{\mathcal{L}}$,
we find that
\[
\int_{U_X}j^3(\phi\circ\phi_X^{-1})^\ast\pounds_{j^3(\xi_\mathcal{C}(\phi
))}\Theta_\mathcal{L} = 0\, .
\]
Then by Theorem \ref{thm_main} and the invariance of $\mathcal{S}$ we have that
\begin{align}\label{101}
0=(\xi_\mathcal{C}\intprod d\mathcal{S})(\phi )&= \int_{\partial U_X}
j^3(\phi\circ\phi_X^{-1})^\ast [j^3(\xi_\mathcal{C}(\phi ))\intprod
\Theta_\mathcal{L}]\nonumber \\ &= -
\int_{U_X}j^3(\phi\circ\phi_X^{-1})^\ast[j^3(\xi_\mathcal{C}(\phi ))\intprod
\Omega_\mathcal{L}]\, .
\end{align}

\begin{defn}
Let $J\in\text{Hom}({\mathfrak g}, T^*{\mathcal C} \otimes \Lambda^n(J^3Y))$
satisfy
\begin{equation}\label{momentum}
j^3(\xi_\mathcal{C}(\phi ))\intprod \Omega_\mathcal{L}= d[J(\xi )(\phi )]
\end{equation}
for all $\xi\in{\mathfrak g}$ and $\phi \in \mathcal{C}$. Then the map
$\mathbb{J}:\mathcal{C}\longrightarrow{\mathfrak g}^\ast$ defined by
\begin{equation}
\langle \mathbb{J}(\phi ), \xi\rangle = J(\xi )(\phi )\, \ \ \forall \ 
\xi \in {\mathfrak g}, \phi\in {\mathcal C},
\end{equation}
is the
{\bf covariant momentum map of the action}.
\end{defn}

With this definition, (\ref{101}) becomes
$\int_{U_X}d[j^3(\phi\circ\phi_X^{-1})^\ast \langle \mathbb{J}(\phi ),
\xi\rangle ] = 0$, and since this holds for any $U_X\subset X$, the integrand
must also vanish; thus,
\begin{equation}\label{noether1}
d[j^3(\phi\circ\phi_X^{-1})^\ast \langle \mathbb{J}(\phi ), \xi\rangle ] =
0\, .
\end{equation}
On the other hand, by Stoke's theorem we may also conclude that
\begin{equation}\label{noether2}
\int_{\partial U_X}j^3(\phi\circ\phi_X^{-1})^\ast \langle \mathbb{J}(\phi ),
\xi\rangle = 0\, .
\end{equation}
Last two statements are equivalent, and we refer to them as the covariant
Noether Theorem.

\section{A multisymplectic approach to the Camassa-Holm equation}
\subsection{The Camassa-Holm equation}

The completely integrable bi-Hamiltonian Camassa-Holm (CH)
equation\footnote{ See \cite{FF} and \cite{CH93}.}
\begin{equation}\label{CH}
u_t-u_{yyt}=-3uu_y+2u_yu_{yy}+uu_{yyy}
\end{equation}
is a model for breaking shallow water waves that admits peaked solitary
traveling waves as solutions (see \cite{CH93}, \cite{CHH94}).  Such
solutions, termed {\it peakons}, develop from any initial data with
sufficiently negative slope, and because of the discontinuities in the
first derivative, these solutions are difficult to numerically 
simulate, particularly in
the case of a {\it peakon--antipeakon collision} (see \cite{CHH94}).

 The multisymplectic framework for the CH
equation is intended to provide a foundation for numerical discretization
schemes that preserve the Hamiltonian structure of this model, even at the
discrete level. After developing the multisymplectic framework for
(\ref{CH}), we shall follow \cite{MPS} and develop the entire discrete
multisymplectic approach to second-order field theories, concentrating on the
discrete CH equation as our model problem. Although we shall only produce the
simplest {\it multisymplectic-momentum} conserving algorithm for this
equation, our construction is completely general and will allow for the
creation of $k$th-order accurate schemes for arbitrarily large $k$.

The CH equation (\ref{CH}) is usually expressed in terms of the Eulerian, 
or spatial
velocity field $u(t,y)$, and is the Euler-Poincar\'{e} equation for the
reduced Lagrangian
\begin{equation}\label{lag}
l(u)=\frac{1}{2}\int(u^2+u_y^2)dy.
\end{equation}

Alternatively, one may express (\ref{CH}) in terms of the Lagrangian
variable $\eta(t,x)$ arising from the solution of
\begin{equation}\label{flow}
\frac{\partial}{\partial t} \eta(x,t) =  u(t,\eta(x,t)).
\end{equation}

The Lagrangian approach to the CH equation is ideally suited to the
multisymplectic variational theory, and we begin by specifying our
fiber bundle $\pi_{XY}:Y\rightarrow X$.
Let $X=S^1\times \mathbb{R}$, and $Y=S^1\times \mathbb{R}\times \mathbb{R}$.
We coordinatize $X$ by $(x^1,x^0)$ (or $(x,t)$) and $Y$ by
$(x^1, x^0, y)$ (or $(x, t, y)$).
A smooth section $\phi\in C^\infty(Y)$ represents a physical field and is
expressed in local coordinates by $(x, t, \eta(x,t))$, where $\eta$ is
the Lagrangian flow solving (\ref{flow}).  The material or Lagrangian velocity
$(\partial / \partial t)\eta (x,t)$ is an element of $T_{\phi(x,t)}Y
= T_{(x,t,y)}Y$, where $y=\eta(x,t)$.

Using (\ref{flow}) together with $u_y = \eta_{tx}/\eta_x$, the Lagrangian
representation for the action may be expressed as
\begin{equation}\label{CHaction}
 \mathcal{S}(\phi )=\int_X\frac{1}{2}(\eta_x\eta_t^2+\eta_x^{-1}
\eta_{tx}^2)dxdt\ .
 \end{equation}
The second jet bundle $J^2Y$ is a nine-dimensional manifold
and $2$-holonomic sections of $J^2Y\longrightarrow X$ have local coordinates
\[
j^2(\phi )=(x, t, \eta (x, t), \eta_x(x, t), \eta_t(x, t), \eta_{xx}
(x, t), \eta_{xt}(x, t), \eta_{tx}(x, t), \eta_{tt}(x, t))\, ,
\]
where for smooth sections $\eta_{xt}(x, t)=\eta_{tx}(x, t)$. The
Lagrangian density $\mathcal{L}:J^2Y\longrightarrow \Lambda^2(X)$ is
expressed as
\[
\mathcal{L}(x^1, x^0, y, y_1, y_0, y_{11}, y_{10}, y_{01}, y_{00})=
L(x^1, x^0, y, y_1, y_0, y_{11}, y_{10}, y_{01}, y_{00})dx^1\wedge dx^0\,
.
\]
For the Camassa-Holm equation the Lagrangian density evaluated along the
second jet of a section $\phi $ is given by
\begin{equation}\label{density}
\mathcal{L}(j^2(\phi ))=\left[\frac{1}{2}(\eta_x\eta_t^2+\eta_x^{-1}
\eta_{tx}^2)\right] dx\wedge dt \, .
\end{equation}

As our Lagrangian (\ref{density}) depends only on  $y_1, y_0$, and
$y_{01}$, the Euler-Lagrange equation (\ref{coordinateEL}) simply becomes
\begin{equation}\label{generall}
-\frac{\partial}{\partial x}\left( \frac{\partial L}{\partial\eta_x}
\right) -\frac{\partial}{\partial t}\left( \frac{\partial L}
{\partial\eta_t}\right)+\frac{\partial^2}{\partial t\partial x}\left(
\frac{\partial L}{\partial\eta_{tx}}\right) =0,
\end{equation}
so that we have the Lagrangian version of the CH equation (\ref{CH}) given by
\begin{equation}\label{ELCH}
\frac{1}{2}\left( \left( \frac{\eta_{tx}}{\eta_x}\right) ^2-
\eta_t^2\right) _x - (\eta_x\eta_t)_t +\left(
\frac{\eta_{tx}}{\eta_x}\right) _{xt} = 0.
\end{equation}
By differentiating $u = (\partial/\partial t)\eta \circ \eta^{-1}$
three times, one may verify that (\ref{ELCH}) is indeed equivalent to
(\ref{CH}).

Now, using (\ref{Cartanagain}) we have that the  Cartan form
$\Theta_\mathcal{L}$ is given by
\begin{align}\label{CHCartan}
\Theta_\mathcal{L} &=\frac{\partial L}{\partial\eta_x}d\eta\wedge dt -
\left( \frac{\partial L}{\partial\eta_t}-D_x\left(
\frac{\partial L}{\partial\eta_{tx}}\right) \right)
d\eta\wedge dx + \frac{\partial L}{\partial\eta_{tx}}d\eta_t\wedge dt
\nonumber \\
&+\left( L-\frac{\partial L}{\partial\eta_x}\eta_x-
\frac{\partial L}{\partial\eta_t}\eta_t-
\frac{\partial L}{\partial\eta_{tx}}\eta_{tx}+D_x\left(
\frac{\partial L}{\partial\eta_{tx}}\right)\eta_t\right)
dx\wedge dt\ ,
\end{align}
or
\begin{align}\label{Cartancontact}
  \Theta_\mathcal{L} &=\frac{\partial L}{\partial\eta_x}
  (d\eta\wedge dt-\eta_xdx\wedge dt)
  +\left( \frac{\partial L}{\partial\eta_t}-D_x\left(
\frac{\partial L}{\partial\eta_{tx}}\right) \right)
(-d\eta\wedge dx-\eta_t dx\wedge dt) \nonumber \\
 &+\frac{\partial L}{\partial\eta_{tx}}(d\eta_t\wedge dt -
 \eta_{tx}dx\wedge dt) + Ldx\wedge dt
\end{align}
if written in terms of the system of contact forms.

\subsection{The multisymplectic form formula for the CH equation}
Marsden, Patrick, and Shkoller in their \cite{MPS} paper have demonstrated
how the multisymplectic form formula for first-order field theories
when
applied to nonlinear wave equations generalizes the notion of
symplecticity given by Bridges in \cite{B97}.  
Using the example of the CH equation, we present below a simple
interpretation of the multisymplectic form
formula for the second-order field theories. We show that the MFF formula
is an intrinsic generalization of the conservation law analogous to the one
in Appendix D of \cite{B97}.

Bridges has introduced the notion of a Hamiltonian system on a
multi-symplectic structure. A multi-symplectic structure $(\mathcal{M},
\omega^1, \ldots , \omega^n, \omega^0)$ consists of a manifold
$\mathcal{M}$, the phase space, and a family of pre-symplectic forms. The
phase space $\mathcal{M}$ is a manifold modeled on $\mathbb{R}^{n+1}$.
A Hamiltonian system on a multi-symplectic
structure is then represented symbolically by $(\mathcal{M}, \omega^1,
\ldots , \omega^n, \omega^0, H)$ with governing equation
\begin{equation}\label{bridges}
\omega^1(\frac{\partial Z}{\partial x^1}, v)+\cdots
+\omega^n(\frac{\partial Z}{\partial x^n}, v)+ \omega^0 (\frac{\partial Z
}{\partial t}, v)=\langle\nabla H(Z), v\rangle
\end{equation}
for all vector fields $v$ on $\mathcal{M}$ where $\langle\cdot,
\cdot\rangle$ is an inner product on $T\mathcal{M}$ and $Z(x^1, \ldots
x^n, t)$ is a curve in $\mathcal{M}$. Bridges has shown that this
formulation is natural for studying wave propagation in open systems.
Bridges, in particular, has obtained the following conservation law in the
case of the wave equation:
\begin{equation}\label{wavelaw}
\frac{\partial }{\partial t}\omega^0(Z_t, Z_x)+ \frac{\partial }{\partial
x}\omega^1(Z_t, Z_x) = 0\ .
\end{equation}
This law generalizes the notion of symplecticity of classical mechanics.

Let us make an appropriate choice of the phase space $\mathcal{M}$ for the
Camassa-Holm equation.  Our choice is entirely governed by the coefficients
in the Cartan form (\ref{CHCartan}).
Since the Lagrangian (\ref{density}) does not explicitly depend on time and
space variables, that is the system is autonomous, we identify sections
$\phi$ of $Y$ with mappings $\eta (x, t)$ from $\mathbb{R}^2$ into
$\mathbb{R}$, and similarly, sections of $J^3Y$ with mappings from
$\mathbb{R}^2$ into $\mathbb{R}^{15}$. The Cartan form (\ref{CHCartan})
suggests to introduce the following momenta:
\begin{equation}\label{Legendre}
\left\{
\begin{array}{ll}
p^x=\frac{\partial L}{\partial\eta_x}\quad   & p^t= \frac{\partial
L}{\partial\eta_t}-D_x\left(\frac{\partial L}{\partial\eta_{tx}} \right)\\
p^{tx}=\frac{\partial L}{\partial\eta_{tx}}\quad   & p^{xx}=p^{tt}
=p^{xt}=0\ .
\end{array}
\right.
\end{equation}
Since $\Theta_\mathcal{L}$ is horizontal over $J^1Y$, the covariant
configuration bundle is really $J^1Y\longrightarrow X$, and one should think
of $(\eta , \eta_x, \eta_t)$ as field variables with each field variable
having conjugate multi-momenta. For example, $p^x, p^t$ function as conjugate
spatial and temporal momenta for the field component $\eta$. Then the
transformation
\[
(\eta , \eta_x, \eta_t, \eta_{xx}, \eta_{xt}, \eta_{tx}, \eta_{tt}, \ldots
)\longmapsto (\eta , \eta_x, \eta_t, p^x, p^t, p^{tx})
\]
defines a mapping from the space of vertical sections of
$J^3Y\longrightarrow X$ into the phase space $\mathcal{M}=\mathbb{R}^6$
modeled over $X=\mathbb{R}^2$. We denote this transformation by
$\mathbb{F}L$. Let us now state the result that connects our paper to
Bridges' theory:

\begin{prop}
The multisymplectic form formula (MFF) yields a multi-symplectic structure
$(\mathcal{M}, \omega^1, \omega^0)$ such that the MFF
formula becomes an
intrinsic generalization of the following conservation law: for
any $V, W$ in $\mathcal{F}$ that are $\pi_{XY}$-vertical
\begin{equation}\label{wavelawgeneral}
\frac{\partial }{\partial x}\omega^1(T\mathbb{F}L\cdot j^3(V),
T\mathbb{F}L\cdot j^3(W))+ \frac{\partial }{\partial
t}\omega^0(T\mathbb{F}L\cdot j^3(V), T\mathbb{F}L\cdot j^3(W))=0\ .
\end{equation}
Moreover, the CH equation in both Eulerian form (\ref{CH}) and
Lagrangian form (\ref{ELCH}) 
is equivalent to the Hamiltonian system of equations on the
multi-symplectic structure with the Hamiltonian defined by
\begin{equation}\label{Hamiltonian}
H=L-p^x\eta_x-p^t\eta_t-p^{tx}\eta_{tx}\ .
\end{equation}
\end{prop}

\begin{proof}
Consider two $\pi_{XY}$-vertical vectors $V$ and $W$ in $\mathcal{F}$. Then
$T\mathbb{F}L\cdot j^3(V)$ and $T\mathbb{F}L\cdot j^3(W)$ are 
vertical-over-$X$ vector fields on $X\times\mathcal{M}\longrightarrow X$ 
whose components we shall denote via
\[
(V^\eta , V^{\eta_x}, V^{\eta_t}, V^{p^x}, V^{p^t}, V^{p^{tx}})\ ,
\]
or just numerate by $(V^1, V^2, V^3, V^4, V^5, V^6)$. Thinking of the
components of the transformation $\mathbb{F}L$ as functions on $J^3Y$, we
immediately see that
\begin{equation}\label{111}
\begin{array}{rl}
V^{p^x}&=dp^x\cdot j^3(V)\equiv j^3(V)[p^x]\\ V^{p^t}&=dp^t\cdot
 j^3(V)\equiv j^3(V)[p^t]\\ V^{p^{tx}}&=dp^{tx}\cdot j^3(V)\equiv
j^3(V)[p^{tx}]
\end{array}
\end{equation}
Using expressions (\ref{Legendre}) we express $\Omega_\mathcal{L}$ as:
\begin{align*}
\Omega_\mathcal{L} & =dp^x\wedge d\eta\wedge dt - dp^t\wedge d\eta\wedge
dx + dp^{tx}\wedge d\eta_t\wedge dt \\
 & - \eta_x dp^x\wedge dx\wedge dt -
\eta_t dp^t\wedge dx\wedge dt - \eta_{tx} dp^{tx}\wedge dx\wedge dt\ .
\end{align*}
Next, combining with (\ref{111}) we obtain that
\begin{align*}
j^3(W)\intprod j^3(V)\intprod\Omega_\mathcal{L}&= \{ V^{p^x}W^\eta -
W^{p^x}V^\eta + V^{p^{tx}}W^{\eta_t} - W^{p^{tx}}V^{\eta_t}\} dt \\
 & - \{V^{p^t}W^\eta - W^{p^t}V^\eta \} dx\ ,
\end{align*}
so that
\begin{align}\label{110}
&\int_{\partial U_X}j^3(\phi\circ\phi_X^{-1})^\ast [j^3(W)\intprod
j^3(V)\intprod\Omega_\mathcal{L}] = \\ &\int_{\partial U_X}(V^4W^1- W^4V^1
+V^6W^3-W^6V^3)dt - (V^5W^1-W^5V^1)dx\ .\nonumber
\end{align}
The integral on the right-hand-side of the above equation leads us to 
introduce two
degenerate skew- symmetric matrices $B_1$, $B_0$ on $\mathbb{R}^6$:
\[
B_1=\left[
\begin{array}{cccccc}
0&0&0&1&0&0\\ 0&0&0&0&0&0\\ 0&0&0&0&0&1\\ -1&0&0&0&0&0\\ 0&0&0&0&0&0\\
0&0&-1&0&0&0
\end{array}\right]
\quad B_0=\left[
\begin{array}{cccccc}
0&0&0&0&1&0\\ 0&0&0&0&0&0\\ 0&0&0&0&0&0\\ 0&0&0&0&0&0\\ -1&0&0&0&0&0\\
0&0&0&0&0&0
\end{array}\right]\ .
\]
To each matrix $B_\nu$, we associate the $2$-form $\omega^\nu$ on
$\mathbb{R}^6$ given by $\omega^\nu(u, v)=\langle B_\nu u, v\rangle\equiv
v^\mathrm{T}B_\nu u$, where $u, v \in \mathbb{R}^6$. With the definition
of $\omega^\nu$ and the use of (\ref{110}), the multisymplectic form formula
(\ref{mff}) becomes, for $U_X\subset X$,
\[
\int_{\partial U_X}\omega^1(T\mathbb{F}L\cdot j^3(V), T\mathbb{F}L\cdot
j^3(W))dt - \omega^0(T\mathbb{F}L\cdot j^3(V), T\mathbb{F}L\cdot
j^3(W))=0\ .
\]
Hence by the Stoke's theorem,
\[
\int_{\partial U_X}\left[ \frac{\partial }{\partial
x}\omega^1(T\mathbb{F}L\cdot j^3(V), T\mathbb{F}L\cdot j^3(W))+
\frac{\partial }{\partial t}\omega^0(T\mathbb{F}L\cdot j^3(V),
T\mathbb{F}L\cdot j^3(W))\right] dx\wedge dt =0\ .
\]
Since $U_X$ is arbitrary, we obtain the desired conservation law
(\ref{wavelawgeneral}).

In the special case, when the components $V^\eta =\eta_x$ and $W^\eta
=\eta_t$, one may verify that
\begin{align*}
T\mathbb{F}L\cdot j^3(V)&=(\eta , \eta_x, \eta_t, p^x, p^t, p^{tx})_{,x}\
,
\\ T\mathbb{F}L\cdot j^3(W)&=(\eta , \eta_x, \eta_t, p^x, p^t,
p^{tx})_{,t}\ ,
\end{align*}
so that, letting $Z$ denote an element $(\eta , \eta_x, \eta_t, p^x, p^t,
p^{tx})\in \mathcal{M}$, the formula (\ref{wavelawgeneral}) takes
the special form 
\[
\frac{\partial }{\partial x}\omega^1(Z_t, Z_x)+ \frac{\partial }{\partial
t}\omega^0(Z_t, Z_x) = 0\ ,
\]
which is the complete analogue of Bridges' conservation law
(\ref{wavelaw}) for the wave equation.

Next, since the inner product $\langle\cdot , \cdot\rangle$ is independent
of $Z\in\mathcal{M}$, the Hamiltonian system of equations (\ref{bridges})
on the multi-symplectic structure $(\mathcal{M}, \omega^1, \omega^0)$ may
be written as
\[
Z_x\intprod\omega^1 + Z_t\intprod\omega^0 = \nabla H\ , \mbox{ or}\quad
B_1Z_x + B_0Z_t = \nabla H\ ,
\]
which results in
\begin{align*}
\frac{\partial }{\partial x}p^x + \frac{\partial }{\partial t}p^t & =
\frac{\partial H}{\partial\eta }\\ \frac{\partial }{\partial x}p^{tx} & =
\frac{\partial H}{\partial\eta_t}\\ \frac{\partial }{\partial x}\eta & = -
\frac{\partial H}{\partial p^x}\\ \frac{\partial }{\partial t}\eta & = -
\frac{\partial H}{\partial p^t}\\ \frac{\partial }{\partial x}\eta_t & = -
\frac{\partial H}{\partial p^{tx}}\ .
\end{align*}
With the choice of the Hamiltonian (\ref{Hamiltonian}) the last four
equations yield identities, and the first equation becomes
\[
\frac{\partial }{\partial x}p^x + \frac{\partial }{\partial t}p^t = 0\ .
\]
Using the Legendre transformation expressions (\ref{Legendre}) for $p^x,
p^t$, the latter equation recovers the Euler-Lagrange equation
(\ref{generall}), hence (\ref{ELCH}). In other words, the Euler-Lagrange
equations on $J^3Y$ are equivalent to Hamilton's equations on the multi-
symplectic structure $(\mathcal{M}, \omega^1, \omega^0, H)$.
\end{proof}

\section{Discrete second-order multisymplectic field theory}
\subsection{A general construction}\label{general}
We shall now generalize the Veselov-type discretization of
first-order field theory given in \cite{MPS} to second-order field theories,
using the Camassa-Holm equation as our example. We discretize $X$ by
$\mathbb{Z} \times\mathbb{Z}=\{(i, j)\}$ and
the fiber bundle $Y$ by $X\times\mathbb{R}$.
Elements of Y over the base point $(i, j)$ are written as $y_{ij}$ and
the projection $\pi_{XY}$ acts on $Y$ by $\pi_{XY}(y_{ij})=(i, j)$. The
fiber over $(i, j)\in X$ is denoted $Y_{ij}$.

For the general case of a second-order Lagrangian one must define the
discrete second jet bundle of $Y$, and this discretization depends on how
one chooses to approximate the partial derivatives of the field.
For example, using central differencing and a fixed timestep $k$ and
spacestep $h$, we have that

\begin{align}\label{center}
  \eta_x & \approx\frac{y_{i+1j}-y_{i-1j}}{2h}\, ,\quad
  \eta_t\approx\frac{y_{ij+1}-y_{ij-1}}{2k}\, ,\quad
  \eta_{xx}\approx\frac{y_{i-1j}-2y_{ij}+y_{i+1j}}{h^2}\, ,\nonumber \\
  \eta_{tx} & \approx\frac{y_{i+1j+1}-y_{i+1j-1}+y_{i-1j-1}-y_{i-1j+1}}
  {4hk}\, ,\quad \eta_{tt}\approx\frac{y_{ij-1}-2y_{ij}+y_{ij+1}}{k^2}\, ,
\end{align}
where $y_{ij}=\eta (x_i, t_j)$ and $\{(x_i, t_j)\}$ form a uniform grid in
continuous spacetime.

\vskip -0.5in
\begin{figure}
\epsfxsize 4.in
\epsfysize 4.in
\hskip 1.25in
\epsfbox{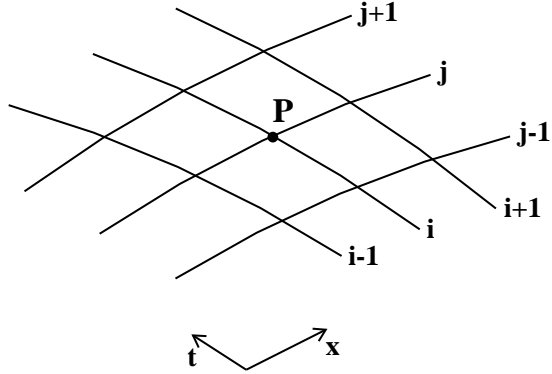}

\vskip -0.5in \caption{Equivalent computational grid in the physical domain.}
\label{FIG41}
\end{figure}

We observe that an $9$-tuple
\[
(y_{i-1j-1},\, y_{i-1j},\, y_{i-1j+1},\, y_{ij-1},\, y_{ij},\, y_{ij+1},\,
 y_{i+1j-1},\,
y_{i+1j},\, y_{i+1j+1})
\]
is sufficient to approximate $j^2\phi (P)$, where $P$ is in the center of
the cell
\begin{eqnarray*}
\boxplus_{ij}\equiv(\begin{array}[t]{ccc}
                    (i-1, j-1),& (i-1, j),& (i-1, j+1),\\
                    (i, j-1),& (i, j),& (i, j+1),\\
                    (i+1, j-1),& (i+1, j),& (i+1, j+1)\ )\, .
                    \end{array}
\end{eqnarray*}
Let
$X^\boxplus $ denote the set of cells, i.e. $X^\boxplus =\{\boxplus_{ij}
\ |\ (i, j)\in X\}$. Components of a cell are called vertices, and are
numbered from first to ninth. A point $(i, j)\in X$ is \textbf{touched} by
a cell if it is a vertex of that cell. If $U\subseteq V$, then $(i, j)\in X$
is an \textbf{interior point} of $U$ if $U$ contains all cells touching $(i, j)$.
The \textbf{interior} $\mathrm{int }\, U$ of $U$ is the set of all interior points of
$U$. The \textbf{closure} $\mathrm{cl }\, U$ of $U$ is the union of all cells
touching interior points of $U$. A \textbf{boundary point} of $U$ is
a point in $U$ and $\mathrm{cl }\, U$ which is not an interior point. The
\textbf{boundary} of $U$ is the set of boundary points, so that
$\partial U\equiv (U\cap \mathrm{cl }\, U)\setminus \mathrm{int }\, U$.

A \textbf{section} of the configuration bundle $Y\longrightarrow X$ is a
map $\phi :U\subseteq X\longrightarrow Y$ such that $\pi_{XY}\circ\phi =
\mathrm{id}_U$. We are now ready to define the discrete
multisymplectic phase space.
\begin{defn}
The discrete \textbf{second jet bundle} of $Y$ is given by
\begin{align*}
J^2Y\equiv \{&(y_{i-1j-1},\, y_{i-1j},\, y_{i-1j+1},\, y_{ij-1},\, y_{ij},\, y_{ij+1},\,
 y_{i+1j-1},\,
y_{i+1j},\, y_{i+1j+1})\ |\\
& (i, j)\in X\, ,\
y_{i-1j-1},\, \ldots ,\, y_{i+1j+1}\in \mathbb{R}\} \\
\equiv\ \  &X^\boxplus \times \mathbb{R}^9\, .
\end{align*}
\end{defn}
The fiber over $(i, j)\in X$ is denoted $J^2Y_{ij}$. We define the
\textbf{second jet extension} of a section $\phi$ to be the map
$j^2\phi :X\longrightarrow J^2Y$ given by
\begin{eqnarray*}
j^2\phi(i, j)\equiv (
\boxplus_{ij},\begin{array}[t]{ccc}
                    \phi (i-1, j-1),& \phi (i-1, j),& \phi (i-1, j+1),\\
                    \phi (i, j-1),& \phi (i, j),& \phi (i, j+1),\\
                    \phi (i+1, j-1),& \phi (i+1, j),& \phi (i+1, j+1)\ )\, .
                    \end{array}
\end{eqnarray*}
Given a vector field $v$ on $Y$ the \textbf{second jet extension} of $v$ is the
vector field $j^2v$ on $J^2Y$ defined by
\begin{eqnarray*}
  j^2v(y_{i-1j-1},\, \ldots ,\, y_{i+1j+1})\equiv
  \begin{array}[t]{ccc}(v(y_{i-1j-1}),&\, v(y_{i-1j}),&\, v(y_{i-1j+1}),\, \\
                     v(y_{ij-1}),&\, v(y_{ij}),&\, v(y_{ij+1}),\, \\
                      v(y_{i+1j-1}),&\, v(y_{i+1j}),&\, v(y_{i+1j+1}))\, .
  \end{array}
\end{eqnarray*}

Of course, this may easily be generalized to more accurate differencing
schemes that require more than nine grid points to define second
partial derivatives.

\subsection{A multisymplectic-momentum algorithm for the CH equation}

Restricting our attention to the CH equation and noting that its
Lagrangian depends only on $\eta_x, \eta_t, \eta_{tx}$, we may
significantly simply our discretization of the second jet bundle $J^2Y$;
this will substantially reduce our calculations and simplify the exposition.

To approximate $j^2\phi (P)$ we choose the forward difference evaluations
of $\eta_x, \eta_t, \eta_{tx}$:
\begin{align*}
  \eta_x &\approx \frac{y_{i+1j}-y_{ij}}{h}\, ,\quad
  \eta_t \approx  \frac{y_{ij+1}-y_{ij}}{k}\, , \\
  \eta_{tx} &\approx \frac{y_{i+1j+1}-y_{i+1j}-y_{ij+1}+y_{ij}}{hk}\, .
\end{align*}

For this particular choice, our cell reduces to a rectangle. A
\textbf{rectangle} $\square$ of $X$ is an ordered $4$-tuple of the
form
\[
\square_{ij}=((i, j),\, (i+1, j),\, (i+1, j+1),\, (i, j+1))\, .
\]
For each rectangle, $\square^1,\, \square^2, \, \square^3, \, \square^4$
stand for the first, second, third, and fourth vertices respectively. If
$(i,j)$ is the first vertex, we shall denote the rectangle by $\square_{ij}$.
The set of all rectangles in $X$ is denoted by $X^\square$. The
set-theoretical definitions of Subsection \ref{general}  apply here. For
example, a point $P=(i, j)\in X$ is \textbf{touched} (see Fig. 4.2) by four
rectangles $\square_{ij},\, \square_{i-1j},\,
 \square_{i-1j-1},\, \square_{ij-1}$, etc.

%\vskip -0.1in

\begin{figure}
\epsfxsize 3.in \epsfysize 3.in \hskip 0.25in \epsfbox{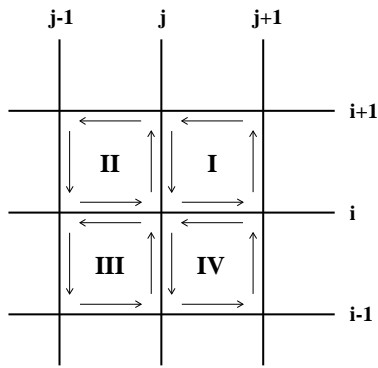}

\vskip -0.5in \caption{The rectangles which touch $(i, j)$.} \label{FIG42}
\end{figure}

%\vskip -0.5in

 Again as (\ref{density}) does not depend on $\eta_{xx},
\eta_{tt}$, we may restrict ourselves to a subbundle $\tilde{\mathcal{B}}$ of
the continuous $J^2Y$ defined via $\tilde{\mathcal{B}}\equiv\{s\in J^2Y\, |\,
s_{\mu\mu }=0\ \ \mbox{for}\ \ \mu = 1, 0\}$. Then the discrete analogue
$\mathcal{B}$ (see Fig. 4.3) of $\tilde{\mathcal{B}}$ is identified with
\begin{equation*}
%\begin{align*}
  \tilde{\mathcal{B}} \equiv\{(y_{ij}, y_{i+1j}, y_{i+1j+1},
  y_{ij+1})\  |\  (i, j)\in X, \  y_{ij}, y_{i+1j}, y_{i+1j+1},
  y_{ij+1}\in \mathbb{R}\}  \equiv X^\square\times\mathbb{R}^4\, .
%  & \equiv X^\square\times\mathbb{R}^4\, .
%\end{align*}
\end{equation*}
\vskip -0.1in
\begin{figure}
\epsfxsize 6.in \epsfysize 4.in \hskip 1.25in \epsfbox{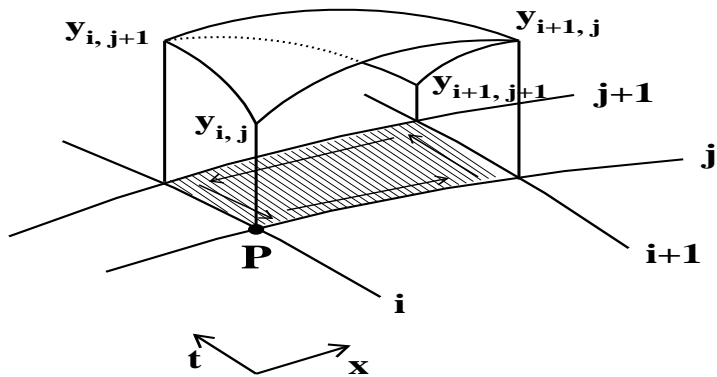}

\vskip -0.5in \caption{Interpretation of an element of $J^2Y$ when $X$ is
discrete. } \label{FIG43}
\end{figure}

 For a section $ \phi :U\subseteq X\longrightarrow Y$, we define
the \textbf{second jet extension} of $\phi $ to $\mathcal{B}$ to be the map
$j^2\phi : U\subseteq X\longrightarrow \mathcal{B}$ via
\[
j^2\phi(i, j)=(\square_{ij},\, \phi (\square^1),\, \phi (\square^2)
, \, \phi (\square^3), \, \phi (\square^4))\, .
\]
Given a vector field $v$ on $Y$ we extend it to a vector field $j^2v$ on
$\mathcal{B}$ by
\[
j^2v(y_{ij}, y_{i+1j}, y_{i+1j+1}, y_{ij+1})=
(v(y_{ij}), v(y_{i+1j}), v(y_{i+1j+1}), v(y_{ij+1}))\, .
\]

A \textbf{discrete Lagrangian} on $\mathcal{B}$ is then a function
$L:\mathcal{B}\longrightarrow \mathbb{R}$ of five variables
$\square_{ij},\, y_1,\, y_2,\, y_3,\, y_4$, where the $y$-variables are
labeled in the order they appear in a $4$-tuple. Let $U$ be a
\textbf{regular} subset of $X$, i.e. $U$ is exactly the union of its
interior and boundary. Let $\mathcal{C}_U$ denote the set of sections
of $Y$ on $U$, so $\mathcal{C}_U$ is the manifold $\mathbb{R}^{|U|}$.
\begin{defn}
The \textbf{discrete action} is a real valued function on $\mathcal{C}_U$
 defined by the rule
 \begin{equation}\label{actionsum}
 \mathcal{S}(\phi )\equiv \sum_{\square\subseteq U;(i, j)=\square^1}
 L\circ j^2\phi (i, j)\ .
 \end{equation}
\end{defn}
Given a section $\phi $ on $U$ acting as $\phi (i, j)=y_{ij}$, one can
define an element $V\in T_\phi \mathcal{C}_U$ to be a map $V:
U\longrightarrow TY$ acting as $V(i, j)= (\phi (i, j),\, v_{ij})$,
where $v_{ij}$ is thought as a vector emanating from $y_{ij}=\phi (i, j)$.
Given an element $V\in T_\phi \mathcal{C}_U$ one can always extend it to a
vector field $v$ on $Y$. On the other hand, given a vector field $v$ on
$Y$ , $V\equiv v\circ\phi $ is an element of $T_\phi \mathcal{C}_U$.
Thus, it is sufficient to work with vector fields $v$ on $Y$ alone.

If $v$ is a vector field on $Y$, consider its restriction
$v\left. \right|_{Y_{ij}}$ to the fiber $Y_{ij}$.  Let $F^v_\lambda :
Y_{ij}\longrightarrow Y_{ij}$ be the flow of $v\left. \right|_{Y_{ij}}$.
Then by definition of the flow,
\[
v(\phi(i, j)) = \left. \frac{d}{d\lambda}\right|
_{\scriptscriptstyle{\lambda =0}}F^v_\lambda (\phi (i, j))\, .
\]
Therefore there is the $1$-parameter family of sections on $U$ defined
by $\phi^\lambda\equiv F^v_\lambda\circ\phi$ such that $\phi^0=\phi $
and $\left. \frac{d}{d\lambda}\right|
_{\scriptscriptstyle{\lambda =0}}\phi^\lambda = v\circ\phi=V$. Thus, the
\textbf{variational principle} is to seek those sections $\phi $ for which
\begin{equation}\label{discvar}
\left. \frac{d}{d\lambda}\right|
_{\scriptscriptstyle{\lambda =0}}\mathcal{S}(F^v_\lambda\circ\phi)=0
\end{equation}
for all vector fields $v$ on $Y$.

\subsection{The discrete Euler-Lagrange equations.}

With our choice of $\mathcal{B}$, the discrete Lagrangian for the
Camassa-Holm equation  is
\begin{equation}\label{DCHL}
L(y_1,\, y_2,\, y_3,\, y_4)=
\frac{1}{2}\left(\frac{y_2-y_1}{h}\cdot\frac{(y_4-y_1)^2}{k^2}
+ \frac{h}{y_2-y_1}\cdot\frac{(y_3-y_2-y_4+y_1)^2}{h^2k^2}\right) \, .
\end{equation}
The variational principle yields the \textbf{discrete Euler-Lagrange
field equations} (DEL equations) as follows. Choose an arbitrary point
$(i, j)\in U$. Henceforth, with a slight abuse of notation, we shall
write $y_{ij}$
for $\phi (i, j)$. The action (\ref{actionsum}), written with its
summands containing $y_{ij}$ explicitly, is (see Figure 3.2)
\begin{align*}
\mathcal{S}=\cdots & + L(y_{ij}, y_{i+1j}, y_{i+1j+1}, y_{ij+1})
+ L(y_{i-1j}, y_{ij}, y_{ij+1}, y_{i-1j+1})\\
 & +  L(y_{i-1j-1}, y_{ij-1}, y_{ij}, y_{i-1j})
 +  L(y_{ij-1}, y_{i+1j-1}, y_{i+1j}, y_{ij}) + \cdots\, .
\end{align*}
Differentiating with respect to $y_{ij}$ yields the DEL equations:
\begin{align*}
& \frac{\partial L}{\partial y_1}(y_{ij}, y_{i+1j}, y_{i+1j+1}, y_{ij+1})
 + \frac{\partial L}{\partial y_2}(y_{i-1j}, y_{ij}, y_{ij+1}, y_{i-1j+1})  \\
& +\frac{\partial L}{\partial y_3}(y_{i-1j-1}, y_{ij-1}, y_{ij}, y_{i-1j})
 +\frac{\partial L}{\partial y_4}(y_{ij-1}, y_{i+1j-1}, y_{i+1j},
 y_{ij})=0
\end{align*}
for all $(i, j)\in \mathrm{int}\ U$. Equivalently, these equations may be
written as
\begin{equation}\label{DELF}
\sum_{l;\square ;(i, j)=\square^l}\frac{\partial L}{\partial y_l}
(\phi (\square^1),\, \phi (\square^2),\, \phi (\square^3),\,
\phi (\square^4))=0
\end{equation}
for all $(i, j)\in \mathrm{int}\ U$. Computing and evaluating
$\frac{\partial L} {\partial y_i}$ along rectangles touching an interior
point $(i, j)$, and substituting these expressions into
(\ref{DELF}), we obtain the discrete Euler-Lagrange equations for the
CH equation:
\begin{align}\label{DELFCH}
&\frac{(\triangle_ky_{i+1j}-\triangle_ky_{ij})^2}{2hk^2(\triangle_hy_{ij})^2}
-\frac{(\triangle_ky_{ij}-\triangle_ky_{i-1j})^2}{2hk^2(\triangle_hy_{i-1j})^2}
-\frac{(\triangle_ky_{ij})^2}{2hk^2}+\frac{(\triangle_ky_{i-1j})^2}{2hk^2}
\nonumber \\
&+\frac{(\triangle_ky_{i+1j}-\triangle_ky_{ij})}{hk^2(\triangle_hy_{ij})}
-\frac{(\triangle_ky_{ij}-\triangle_ky_{i-1j})}{hk^2(\triangle_hy_{i-1j})}
-\frac{(\triangle_ky_{i+1j-1}-\triangle_ky_{ij-1})}{hk^2(\triangle_hy_{ij-1})}
 \\
&+\frac{(\triangle_ky_{ij-1}-\triangle_ky_{i-1j-1})}{hk^2(\triangle_hy_{i-1j-1})}
-\frac{(\triangle_hy_{ij})(\triangle_ky_{ij})}{hk^2}
+\frac{(\triangle_hy_{ij-1})(\triangle_ky_{ij-1})}{hk^2}=0\, ,\nonumber
\end{align}
where
\begin{align*}
 \triangle_ky_{ij}&= y_{ij+1}-y_{ij}, \quad \quad\triangle_ky_{i-1j}
 = y_{i-1j+1}-y_{i-1j}\, ,\\
\triangle_hy_{ij}&= y_{i+1j}-y_{ij},\quad \quad\triangle_ky_{i+1j}
 = y_{i+1j+1}-y_{i+1j}\, .
\end{align*}

To see that (\ref{DELFCH}) are indeed approximating the continuous
Euler-Lagrange equation (\ref{ELCH}), notice that the first two
terms combine to approximate
\begin{align*}
\frac{1}{2}\left(\left(\frac{\eta_{tx}}{\eta_x}\right)^2\right)_x
\approx
\frac{1}{2}\ \frac{1}{h}
&\left[
\left(\frac{\triangle_ky_{i+1j}-
\triangle_ky_{ij}}{hk} \right)^2\left/\left(
\frac{\triangle_hy_{ij}}{h}\right)\right.^2\right. \\
& - \left.\left(\frac{\triangle_ky_{ij}-
\triangle_ky_{i-1j}}{hk} \right)^2\left/\left(
\frac{\triangle_hy_{i-1j}}{h}\right)\right.^2
\right]\, .
\end{align*}
As to the third and fourth terms of (\ref{DELFCH}),
\[
-\frac{1}{2}(\eta_t^2)_x\approx - \frac{1}{2}\ \frac{1}{h}
\left[
\left(\frac{\triangle_ky_{ij}}{k}\right)^2
-\left(\frac{\triangle_ky_{i-1j}}{k}\right)^2
\right]\, .
\]
Next, the fifth, sixth, seventh, and eighth terms combine as
\begin{align*}
\left(\frac{\eta_{tx}}{\eta_x}\right)_{tx}\approx
\frac{1}{hk}\left[
\left(\frac{\triangle_ky_{i+1j}-
\triangle_ky_{ij}}{hk}\left/
\frac{\triangle_hy_{ij}}{h}\right.\right)
-\left(\frac{\triangle_ky_{ij}-
\triangle_ky_{i-1j}}{hk}\left/
\frac{\triangle_hy_{i-1j}}{h}\right.\right)
\right. \\
-\left.
\left(\frac{\triangle_ky_{i+1j-1}-
\triangle_ky_{ij-1}}{hk}\left/
\frac{\triangle_hy_{ij-1}}{h}\right.\right)
+\left(\frac{\triangle_ky_{ij-1}-
\triangle_ky_{i-1j-1}}{hk}\left/
\frac{\triangle_hy_{i-1j-1}}{h}\right.\right)
\right]\, .
\end{align*}
Finally, the last two terms of (\ref{DELFCH}) approximate
\[
-(\eta_x\eta_t)_t\approx -\frac{1}{k}
\left(
\frac{\triangle_hy_{ij}}{h}\cdot\frac{\triangle_ky_{ij}}{k}
-\frac{\triangle_hy_{ij-1}}{h}\cdot\frac{\triangle_ky_{ij-1}}{k}
\right)\, .
\]

The numerical scheme (~\ref{DELFCH}) proceeds as follows:
Suppose that
\[
\triangle_hy_{ij},\ \triangle_hy_{i-1j},\ \triangle_hy_{i-1j-1},\
\triangle_hy_{ij-1},\  \triangle_ky_{ij-1},\  \triangle_ky_{i-1j-1},\
\triangle_ky_{i+1j-1}
\]
as known from  the two previous time steps; then (\ref{DELFCH})
may be written as
\[
\mathbf{F}(\triangle_ky_{ij},\ \triangle_ky_{i+1j},\ \triangle_ky_{i-1j})
=0\, .
\]
These are implicit equations which must be solved for
$y_{ij+1},\ 1\leq i\leq N$, where $N$ is the size of the spatial grid.

\subsection{The discrete Cartan form.}

We consider arbitrary variations which are in no way constrained
on the boundary $\partial U$. For each $(i, j)\in
\partial U$ there is at least one rectangle in $U$ touching $(i, j)$
since $(i, j)\in \mathrm{cl}\ U$ and $U$ is regular. On the other
hand, not all four rectangles touching $(i, j)$ are in $U$ since
$(i, j)\not\in\mathrm{int}\ U$. Therefore, each $(i, j)\in\partial U$
occurs as the $l^{\mathrm{th}}$ vertex for either one, two, or three of
the $l\in {1, 2, 3, 4}$ and the corresponding $l^{\mathrm{th}}$ boundary
expressions are given by
\begin{align}\label{boundarylist}
&\frac{\partial L}{\partial y_1}(y_{ij}, y_{i+1j}, y_{i+1j+1},
y_{ij+1})  V(i, j)\, ,\nonumber\\
 & \frac{\partial L}{\partial y_2}(y_{i-1j}, y_{ij}, y_{ij+1},
 y_{i-1j+1}) V(i, j) \, ,\nonumber \\
&  \frac{\partial L}{\partial y_3}(y_{i-1j-1}, y_{ij-1}, y_{ij},
y_{i-1j}) V(i, j)\, ,\nonumber \\
& \frac{\partial L}{\partial y_4}(y_{ij-1}, y_{i+1j-1}, y_{i+1j},
 y_{ij})V(i, j)\, ,
\end{align}
where $y_{ij}=\phi (i, j)$.
The sum of all such terms is the contribution to $d\mathcal{S}$ from
the boundary $\partial U$.
We thus define the \textit{four} $1$-forms on $\mathcal{B}
\subseteq J^2Y$ by
\begin{align*}
\Theta_L^1 (y_{ij}, y_{i+1j}, y_{i+1j+1}, y_{ij+1})&\cdot
(v_{y_{ij}}, v_{y_{i+1j}}, v_{y_{i+1j+1}}, v_{y_{ij+1}}) \\ &\equiv
  \frac{\partial L}{\partial y_1}(y_{ij}, y_{i+1j}, y_{i+1j+1},
y_{ij+1})\cdot (v_{y_{ij}},\ 0,\ 0,\ 0)\, ,\\
\Theta_L^2 (y_{ij}, y_{i+1j}, y_{i+1j+1}, y_{ij+1})&\cdot
(v_{y_{ij}}, v_{y_{i+1j}}, v_{y_{i+1j+1}}, v_{y_{ij+1}}) \\ &\equiv
  \frac{\partial L}{\partial y_2}(y_{ij}, y_{i+1j}, y_{i+1j+1},
y_{ij+1})\cdot (0,\ v_{y_{i+1j}},\ 0,\ 0)\, ,\\
\Theta_L^3 (y_{ij}, y_{i+1j}, y_{i+1j+1}, y_{ij+1})&\cdot
(v_{y_{ij}}, v_{y_{i+1j}}, v_{y_{i+1j+1}}, v_{y_{ij+1}}) \\ &\equiv
  \frac{\partial L}{\partial y_3}(y_{ij}, y_{i+1j}, y_{i+1j+1},
y_{ij+1})\cdot (0,\ 0,\ v_{y_{i+1j+1}},\ 0)\, ,\\
\Theta_L^4 (y_{ij}, y_{i+1j}, y_{i+1j+1}, y_{ij+1})&\cdot
(v_{y_{ij}}, v_{y_{i+1j}}, v_{y_{i+1j+1}}, v_{y_{ij+1}}) \\ &\equiv
  \frac{\partial L}{\partial y_4}(y_{ij}, y_{i+1j}, y_{i+1j+1},
y_{ij+1})\cdot (0,\ 0,\ 0,\ v_{y_{ij+1}})\, .
\end{align*}
We regard the $4$-tuple $(\Theta_L^1, \Theta_L^2, \Theta_L^3,
\Theta_L^4)$ as being the discrete analogue of the multisymplectic
form $\Theta_\mathcal{L}$.
Given a vector field $v$ on $Y$ such that $V=v\circ\phi $,
the first expression from the list (\ref{boundarylist}) becomes
$[(j^2\phi )^\ast (j^2v\intprod\Theta_L^1)]((i, j))$, the others
written similarly. With this notation, $d\mathcal{S}$ may be expressed
 as
\begin{align}\label{ddS}
d\mathcal{S}(\phi )\cdot V &=
\sum_{(i, j)\in\mathrm{int}\ U}
\left(
\sum_{\stackrel{\square\subseteq U; l;}{(i, j)=\square^l}}
[(j^2\phi )^\ast(j^2v\intprod\Theta_L^l)](\square^1)
\right)   \nonumber \\
&+\sum_{(i, j)\in\partial U}
\left(
\sum_{\stackrel{\square\subseteq U; l;}{(i, j)=\square^l}}
[(j^2\phi )^\ast(j^2v\intprod\Theta_L^l)](\square^1)
\right)\, .
\end{align}

\subsection{The discrete multisymplectic form formula.}
For a rectangle $\square$ in $X$, define the projection $\pi_\square :
\mathcal{C}_U\longrightarrow\mathcal{B}$ by
\[
\pi_\square (\phi )\equiv (\square , \phi (\square^1),
\phi (\square^2), \phi (\square^3), \phi (\square^4))\, .
\]
Calculating the form $\pi_\square^\ast\Theta_L^l$ on
$\mathcal{C}_U$ gives
\[
(\pi_\square^\ast\Theta_L^l)(\phi )\cdot V =
\frac{\partial L}{\partial y_l}
(\phi (\square^1),\, \phi (\square^2),\, \phi (\square^3),\,
\phi (\square^4))\ V(\square^l)\, .
\]
This immediately implies that the variation (\ref{ddS}) can be written as
\begin{align}\label{dSS}
d\mathcal{S}(\phi )\cdot V &=
\sum_{(i, j)\in\mathrm{int}\ U}
\left(
\sum_{\stackrel{\square\subseteq U; l;}{(i, j)=\square^l}}
(\pi_\square^\ast\Theta_L^l)(\phi )\cdot V
\right) \nonumber  \\
&+\sum_{(i, j)\in\partial U}
\left(
\sum_{\stackrel{\square\subseteq U; l;}{(i, j)=\square^l}}
(\pi_\square^\ast\Theta_L^l)(\phi )\cdot V
\right)\, .
\end{align}
Define the $1$-forms $\alpha_1$ and $\alpha_2$ on the space
of sections $\mathcal{C}_U$ to be the first and the second terms
on the right hand side
of (\ref{dSS}), respectively.

As in the previous section we would like to derive the discrete
analogue of symplecticity of the flow in mechanics.
Let $\phi^\lambda$ be a curve of solutions of (~\ref{DELF}) that
passes through $\phi$ at zero with $V=\left.\frac{d}{d\lambda }
\right| _{\scriptscriptstyle{\lambda =0}}\phi^\lambda$. Then
for each interior point $(i, j)$ and each $\lambda$, the following
holds:
\[
\sum_{l;\square ;(i, j)=\square^l}\frac{\partial L}{\partial y_l}
(\phi^\lambda (\square^1),\, \phi^\lambda (\square^2),\, \phi
^\lambda (\square^3),\, \phi^\lambda (\square^4))=0\, .
\]
Differentiating these equations with respect to $\lambda$ at
$\lambda =0$, we obtain the
\begin{defn}
If $\phi$ is a solution of the discrete Euler-Lagrange
equations (~\ref{DELF}), then a
\textbf{first-variation} equation solution at $\phi$ is
a vector $V\in T_\phi\mathcal{C}_U$ such that for each
$(i, j)\in\mathrm{int}\ U$,
\begin{equation}\label{varfirst}
\sum_{l;\square ;(i, j)=\square^l}
\sum_{k=1}^4
\frac{\partial^2 L}{\partial y_k \partial y_l}
(\phi (\square^1),\, \phi (\square^2),\, \phi (\square^3),\,
\phi (\square^4))\ V(\square^k)=0\, .
\end{equation}
\end{defn}

By definition of the forms $\alpha_1$ and $\alpha_2$,
$d\mathcal{S}=\alpha_1+\alpha_2$. Since $d^2\mathcal{S}=0$,
$d\alpha_1+d\alpha_2=0$.
Using (\ref{dSS}) and denoting the vertices of $\square$
by $y_1, y_2, y_3, y_4$, we have that
$\Theta_L^l=\frac{\partial L}{\partial y_l}dy_l$, which
implies that for all $l=1, 2, 3, 4$
\[
\Omega_L^l=\sum_{k=1}^4\frac{\partial^2L}{\partial y_k\partial y_l}
dy_k\wedge dy_l\, .
\]
Therefore,
\begin{align}\label{pullbackomega}
\pi_\square^\ast\Omega_L^l(\phi ) (V, W)
&= \Omega_L^l(\pi_\square(\phi )) (T_\phi\pi_\square\cdot V,
T_\phi\pi_\square\cdot W)\nonumber \\
&= \Omega_L^l(\phi (\square^1),\, \cdots ,\,
\phi (\square^4))\cdot((V(\square^1),\, \cdots ,\, V(\square^4)),\
(W(\square^1),\, \cdots ,\,W(\square^4)))\nonumber \\
&=\sum_{k=1}^4\frac{\partial^2L}{\partial y_k\partial y_l}
(\phi (\square^1),\, \phi (\square^2),\, \phi (\square^3),\,
\phi (\square^4))\{ V(\square^k)W(\square^l)-V(\square^l)W(\square^k)\}
\, .
\end{align}
Substitution of (\ref{pullbackomega}) into the exterior derivative
of the right hand side of (\ref{dSS}) yields
\begin{align*}
&d\alpha_1(\phi )(V, W)=\sum_{(i, j)\in\mathrm{int}\ U}
\left(
\sum_{\stackrel{\square\subseteq U; l;}{(i, j)=\square^l}}
\sum_{k=1}^4\frac{\partial^2L}{\partial y_k\partial y_l}
\Big( \phi (\square^1),\, \cdots ,\, \phi (\square^4)\Big)
\Big( V(\square^k)W(\square^l)-V(\square^l)W(\square^k)\Big)
\right) \\
&\mbox{and} \\
&d\alpha_2(\phi )(V, W)=\sum_{(i, j)\in\partial U}
\left(
\sum_{\stackrel{\square\subseteq U; l;}{(i, j)=\square^l}}
\sum_{k=1}^4\frac{\partial^2L}{\partial y_k\partial y_l}
\Big( \phi (\square^1),\, \cdots ,\, \phi (\square^4)\Big)
\Big( V(\square^k)W(\square^l)-V(\square^l)W(\square^k)\Big)
\right)\, .
\end{align*}
When specialized to two first-variation solutions $V$ and $W$ at $\phi $,
$d\alpha_1(\phi )(V, W)$ vanishes, because for each interior point
$(i, j)$ all four rectangles
touching it are contained in $U$, and $V(\square^l)=V(i, j)$ and
$W(\square^l)=W(i, j)$. Therefore, $d\alpha_1 = 0$ and the equation
$d^2\mathcal{S}=0$ becomes
$d\alpha_2=0$, which in turn is equivalent to
\begin{equation}\label{DMFF}
\sum_{(i, j)\in\partial U}
\left(
\sum_{\stackrel{\square\subseteq U; l;}{(i, j)=\square^l}}
[(j^2\phi)^\ast(j^2w\intprod j^2v\intprod \Omega_L^l)](\square^1)
\right) =0
\end{equation}
for all vector fields $v, w$ on $Y$.
This is the discrete analogue of the multisymplectic form formula for
the continuous spacetime.

We observe that  $dL=\Theta_L^1+\Theta_L^2+\Theta_l^3+\Theta_l^4$,
which shows that
\[
\Omega_L^1+\Omega_L^2+\Omega_L^3+\Omega_L^4=0\, ,
\]
which in turn implies that only three of the $2$-forms $\Omega_L^l$,
$l=1, 2, 3, 4$ are in fact independent.
In addition, this implies that for a \textit{given and fixed}
rectangle $\square $,
\begin{align*}
0&=\sum_{l=1}^4\pi_\square^\ast\Omega_L^l\ (\phi )(V, W) \\
&=\sum_{l=1}^4\sum_{k=1}^4\frac{\partial^2L}{\partial y_k\partial y_l}
\Big(\phi (\square^1),\, \cdots ,\, \phi (\square^4)\Big)
\Big( V(\square^k)W(\square^l)-V(\square^l)W(\square^k)\Big) \,
\end{align*}
for all sections $\phi$ and all vectors $V, W$.

\subsection{The discrete Noether theorem.}
We would like to derive the discrete version of the Noether theorem
for second order field theories. This is not the most general form
possible as we are working with a particular example. However it
is such as to facilitate the derivation of any other case without
significant effort.

Suppose that a Lie group $G$ with a Lie algebra ${\mathfrak g}$ acts on $Y$
by vertical symmetries such that the Lagrangian $L$ is invariant under the
action. Vertical action simply means that the base elements from $X$ are not
altered under the action, hence the action restricts to each fiber of $Y$.
Let $\Phi :G\times Y\longrightarrow Y$ denote the action of $G$ on $Y$. For
every $g\in G$ let $\Phi_g:Y\longrightarrow Y$ be given by $y_{ij}\longmapsto
\Phi (g, y_{ij})$. We also use the notation $g\cdot y=\Phi_g(y)$ for the
action. Then there is an induced action of $G$ on $\mathcal{B}$ defined in a
natural way:
\[
g\cdot (y_1, y_2, y_3, y_4)=(\Phi (g, y_1), \Phi (g, y_2), \Phi (g,
y_3), \Phi (g, y_4))\, .
\]
Recall that the infinitesimal generator of an action (of a Lie group $G$ on a
manifold $M$) corresponding to a Lie algebra element $\xi\in{\mathfrak g}$ is
the vector field $\xi_M$ on $M$ obtained by differentiating the action with
respect to $g$ at the identity in the direction $\xi$. By the chain rule,
\[
\xi_M(z)=\left.\frac{d}{dt}\right|_{\scriptscriptstyle{t=0}}
[\mathrm{exp}(t\xi)\cdot z]\, ,
\]
where $\mathrm{exp}$ is the Lie algebra exponential map.

Using this formula, we immediately see that
\[
\xi_\mathcal{B}(y_1, y_2, y_3, y_4)=
(\xi_Y(y_1), \xi_Y(y_2), \xi_Y(y_3), \xi_Y(y_4))\, .
\]
The invariance of the Lagrangian under the action implies that
\[
\xi_\mathcal{B}\intprod dL=0\quad \mbox{for all} \ \ \xi\in {\mathfrak g}\, ,
\]
which, for a given $\square$, is equivalent to
\begin{equation}\label{Linvariance}
\sum_{l=1}^4 \frac{\partial L}{\partial y_l}(y_1, y_2, y_3, y_4)
\xi_Y(y_l)=0
\end{equation}
for all $\xi\in {\mathfrak g}$ and all $(y_1, y_2, y_3, y_4)\in \mathcal{B}$.
For each $l$, let us denote by $\pi_{\square^l}: \mathcal{B}\longrightarrow
Y$ the projection onto the $l^{\mathrm{th}}$
 component. Using this projection the four components of the
infinitesimal generator $\xi_{\mathcal{B}}$ are expressed as
\[
\xi_{\mathcal{B}}=\sum_{l=1}^4\xi_{\mathcal{B}}^l=\sum_{l=1}^4
(\xi_Y\circ\pi_{\square^l})\frac{\partial }{\partial y_l}\, .
\]
Hence, the equation (~\ref{Linvariance}) becomes
\begin{equation}\label{Linvariance2}
\sum_{l=1}^4\xi_{\mathcal{B}}^l\intprod\Theta_L^l =0\quad \mbox{for all} \ \
\xi\in{\mathfrak g}\, .
\end{equation}
We observe that for each $l$,
\[
\xi_{\mathcal{B}}^l\intprod\Theta_L^l
=\frac{\partial L}{\partial y_l}\cdot
(\xi_Y\circ\pi_{\square^l})
\]
is a function on $\mathcal{B}$ which we denote by
$J^l(\xi)$. Notice that
$J^l(\xi)=\xi_{\mathcal{B}}^l\intprod\Theta_L^l$
is the discrete multisymplectic analogue of $\xi_M\intprod\omega_L =d
J(\xi )$ in classical mechanics so that $\xi_M$ is the global
Hamiltonian vector field of $J(\xi )$. Many symmetry groups act by
special canonical transformations, i.e. $\pounds_{\xi_M}\theta_L=0$,
in which case $J(\xi )=\xi_M\intprod\theta_L$. In such case, $J(\xi )$
is uniquely defined.

Since $\xi_\mathcal{B}$ is linear in $\xi$, so are the functions
$J^l(\xi )$, and we can replace the Lie group action by a Lie
algebra action $\xi\longmapsto\xi_\mathcal{B}$. Finally, we are ready
to define the momentum maps.
\begin{defn}
There are four ${\mathfrak g}^\ast$-valued \textbf{momentum mappings}
$\mathbb{J}^l,\ l=1, 2, 3, 4$ on $\mathcal{B}$ defined by
\begin{equation}\label{momentummap}
\langle
\mathbb{J}^l(y_1, y_2, y_3, y_4), \xi
\rangle
= J^l(\xi )(y_1, y_2, y_3, y_4)
\end{equation}
for all $\xi\in{\mathfrak g}$ and $(y_1, y_2, y_3, y_4) \in\mathcal{B}$,
where $\langle\cdot , \cdot\rangle$ is the duality pairing.
\end{defn}
The equation (~\ref{Linvariance2}) implies that
\[
\mathbb{J}^1+\mathbb{J}^2+\mathbb{J}^3+\mathbb{J}^4=0\, ,
\]
so, as in the case of the Lagrangian $2$-forms, only three
of the four momenta are essentially distinct.

The discrete version of the Noether theorem for second-order field
theories now follows.
Define the action of the Lie group $G$ on $\mathcal{C}_U$ by
\[
g\cdot\phi\equiv\Phi_g\circ\phi\, ,\ \mbox{i.e.}\quad
(g\cdot\phi)(i, j)=\Phi(g, \phi(i, j))\, ;
\]
since the Lagrangian is $G$-invariant, then
\begin{align*}
\mathcal{S}(g\cdot\phi)&=\sum_{\square\subseteq U}
 L\circ j^2(g\cdot\phi ) (\square^1)\\
&=\sum_{\square\subseteq U}L(g\cdot\phi (\square^1), \cdots ,
g\cdot\phi (\square^4))\\
&=\sum_{\square\subseteq U}L(\phi (\square^1), \cdots ,
\phi (\square^4))=\mathcal{S}(\phi )\, .
\end{align*}
Once again letting $g=\mathrm{exp}(t\xi)$ and differentiating with
respect to $t$ at $t=0$, we obtain that $(\xi_{\mathcal{C}_U}
\intprod d\mathcal{S})(\phi )=0$ for all $\phi\in\mathcal{C}_U$.
One can readily verify that $\xi_{\mathcal{C}_U}(\phi )=
\xi_Y\circ\phi$, which is an element of $T_\phi\mathcal{C}_U$. Thus,
\begin{equation}\label{Sinvariance}
d\mathcal{S}(\phi )\cdot (\xi_Y\circ\phi)=0
\end{equation}
for all $\xi\in{\mathfrak g}$ and $\phi\in\mathcal{C}_U$. Since $\mathcal{S}$
is $G$-invariant, then $G$ sends critical points of $\mathcal{S}$ to
themselves, or in other words, the action restricts to the space of solutions
of the Euler-Lagrange equations. Therefore, if $\phi$ is a solution, so is
$\phi^t\equiv\mathrm{exp}(t\xi )\cdot \phi $, where $\phi^0=\phi$ and
$\left.\frac{d}{dt}\right|_ {\scriptscriptstyle{t=0}}\phi^t=\xi_Y\circ\phi$.
Substituting $\phi^t$ into the discrete Euler-Lagrange equations and
differentiating with respect to $t$ at $t=0$, we obtain that for any $\xi$
and $\phi$, $\xi_Y\circ\phi$ is a first-variation equation solution. Using
(\ref{ddS}), (\ref{Sinvariance}) becomes

\begin{align*}
0=d\mathcal{S}(\phi )\cdot (\xi_Y\circ\phi) &=
\sum_{(i, j)\in\partial U}
\left(
\sum_{\stackrel{\square\subseteq U; l;}{(i, j)=\square^l}}
\frac{\partial L}{\partial y_l}
(\phi (\square^1),\, \cdots ,\, \phi (\square^4))\ \xi_Y
\circ\phi (\square^l)
\right)\nonumber \\
&=\sum_{(i, j)\in\partial U}
\left(
\sum_{\stackrel{\square\subseteq U; l;}{(i, j)=\square^l}}
(\xi_Y\intprod\Theta_L^l)
(\phi (\square^1),\, \cdots ,\, \phi (\square^4))
\right)\nonumber \\
&=\sum_{(i, j)\in\partial U}
\left(
\sum_{\stackrel{\square\subseteq U; l;}{(i, j)=\square^l}}
\mathbb{J}^l
(\phi (\square^1),\, \cdots ,\, \phi (\square^4))(\xi )
\right)
\end{align*}
for all $\phi $ from the solution space and all $\xi $.
Thus, the discrete version of the Noether theorem is
\begin{equation}\label{Noether}
\sum_{(i, j)\in\partial U}
\left(
\sum_{\stackrel{\square\subseteq U; l;}{(i, j)=\square^l}}
[(j^2\phi )^\ast\mathbb{J}^l](\square^1)
\right) =0
\end{equation}
for all $\phi$ from the solution space.

\section*{Acknowledgments}
The authors would like to thank the Center for Nonlinear Studies at
Los Alamos, where most of this work was completed, for providing a wonderful
working environment.  SS was partially supported by NSF-KDI grant 
ATM-98-73133 and the DOE.

\end{document}